\documentclass{conm-p-l}

\def\fs{\footnotesize}          \def\sms{\small}
            \def\ls{\large}

\newcommand{\ad}{\mathop{\rm ad}\nolimits}

\newtheorem{theorem}{Theorem}[section]

\newtheorem{proposition}[theorem]{Proposition}

\theoremstyle{definition}
\newtheorem{definition}[theorem]{definition}

\theoremstyle{remark}
\newtheorem{remark}[theorem]{Remark}

\numberwithin{equation}{section}

\begin{document}

\title[From Quantum Affine Kac-Moody Algebras to Yangians]
{From Quantum Affine Kac-Moody Algebras to Drinfeldians and Yangians}
\author{V.N. Tolstoy}
\address{Institute of Nuclear Physics, Moscow State University,
Moscow 119992, Russia}
\email{tolstoy@nucl-th.sinp.msu.ru}
\thanks{This work was supported by Russian Foundation for Fundamental Research,
grant No. RFBR-02-01-00668 and INTAS OPEN 00-00055.}



\begin{abstract}
A general scheme of construction of Drinfeldians and Yangians 
from quantum non-twisted affine Kac-Moody algebras is presented.
Explicit description of Drinfeldians and Yangians for all Lie algebras of the
classical series $A$, $B$, $C$, $D$ is given in terms of a Chevalley basis.
\end{abstract}

\maketitle
\section{Introduction}
Rational and trigonometric deformations of a universal enveloping $\,$algebra
$U(\mathfrak{g}[u])$ ($\mathfrak{g}$ is a finite-dimensional complex simple
Lie algebra) - Yangian $Y_{\eta}(\mathfrak{g})$ and a quantum ($q$-deformed)
affine algebra $U_q(\mathfrak{g}[u])$ - are playing 
increasingly a role in the theory of integrable systems and the quantum field
theory. For applications it is useful to know various realizations (or
bases) of these deformations. Among various realizations the Chevalley
basis plays a special role due to its compactness and simplicity.
In the case of the quantum affine algebras the Chevalley basis was
originally introduced in the defining relations while the Yangians
were initially defined in terms of a finite Cartesian basis (the
first Drinfeld realization \cite{D1}) as well as in terms of a special
infinite basis (the second Drinfeld realization \cite{D2}). It was
not clear whether does exist a Chevalley basis for the Yangians.

Recently in \cite{T2, T3} it was shown that the Yangian can be
realized in terms of Chevalley basis and, moreover, there was found
a simple connection of this realization with the corresponding
quantum affine algebra (also see \cite{T4, T5}). This connection
allows to introduce a new Hopf algebra called Drinfeldian which
depends on two deformation parameters $q$ and $\eta$. The Drinfeldian
$D_{q\eta}(\mathfrak{g})$ can be considered as quantization of
$U(g[u])$ in the direction of a classical $r$-matrix which is a sum
of simplest rational and trigonometric $r$-matrices. When the parameter
$q$ goes to 1 we obtain the Yangian $Y_{\eta}(\mathfrak{g})$ in terms
of the Chevalley basis.

In this paper we remind once more the procedure of construction of
the Drinfeldian $D_{q\eta}(\mathfrak{g})$ and the Yangian
$Y_{\eta}(\mathfrak{g})$ in terms of the Chevalley basis starting from
the defining relations of the quantum affine algebra
$U_q(\mathfrak{g}[u])$. General defining relations
for $D_{q\eta}(\mathfrak{g})$ depend on some vector
$\tilde{e}_{-\theta}\in U_q(\mathfrak{g})$, where $\theta$ is a maximal
root of the Lie algebra $\mathfrak{g}$. We can choose such vector
$\tilde{e}_{-\theta}$ in order to obtain more simple explicit formulas
for the defining relations of $D_{q\eta}(\mathfrak{g})$ and
$Y_{\eta}(\mathfrak{g})$. Here we give these explicit defining relations
of Drinfeldians and Yangians for all Lie algebras of the classical series
$A$, $B$, $C$, $D$.

The plan of this paper is as follows. In Section 2 a fundamental
result of the paper by A.~Belavin and V.~Drinfeld \cite{BD} concerning
classification of the classical Yang-Baxter equation (CYBE) is briefly
presented. In Section 3 we demonstrate that a sum of the simplest
rational and trigonometric $r$-matrices is also a solution of CYBE.
In order to compare our realization of the Yangian with Drinfeld's
realizations  we present two Drinfeld's realizations in explicit form
in Section 4. Section 5 deals with the quantum algebra
$U_q(\mathfrak{g}[u])$ for an arbitrary simple finite-dimensional complex
Lie algebra $\mathfrak{g}$.
In Section 6 we show how to construct a two-parameter deformation or
the Drinfeldian $D_{q\eta}(\mathfrak{g})$ starting from $U_q(\mathfrak{g}[u])$.
At $q=1$ the Drinfeldian $D_{q\eta}(\mathfrak{g})$ coincides with the Yangian
$Y_{\eta}(\mathfrak{g})$. In Sections 7--10 we give the explicit description
of Drinfeldians and Yangians in the terms of the Chevalley basis for
all Lie algebras of the classical series $A$, $B$, $C$ and $D$.

\setcounter{equation}{0}
\section{Classical Yang-Baxter equation (CYBE).\\
Belavin-Drinfeld classification of CYBE solutions}
The classical Yang-Baxter equation (CYBE) is a equation of the form
\begin{equation}
[r^{12}(u_1,u_2),\,r^{13}(u_1,u_3)+r^{23}(u_2,u_3)]+
[r^{13}(u_1,u_3),\,r^{23}(u_2,u_3)]=0~,
\label{CYBE1}
\end{equation}
where $r^{}(u,v)$ is a function of two variables $u$, $v$ with values in
$\mathfrak{g}\otimes \mathfrak{g}$ ($\mathfrak{g}$ is a finite-dimensional
complex simple Lie algebra). Let $\{e_i\}$ be a basis of $\mathfrak{g}$ then
$r^{}(u,v)$ has the form
\begin{equation}
r(u,v)=\sum_{i,j}f_{ij}(u,v)\,e_i\otimes e_j~,
\label{CYBE2}
\end{equation}
where $f_{\!ij}(u,v)$ are usual functions of the variables $u$ and $v$.
The notations $r^{12}(u_1\!,u_2)$, $r^{13}(u_1,u_3)$ and
$r^{23}(u_2,u_3)$ are defined as follows:
\begin{equation}
r^{12}(u_1,u_2)=\sum_{i,j}f_{ij}(u_1,u_2)\,e_i\otimes e_j\otimes1 
\label{CYBE3}
\end{equation}
etc. 
One usually assumes that the function $r(u,v)$ satisfies the
following additional properties:
\begin{itemize}
\item[{\it(i)}]
the unitary condition: $r^{12}(u,v)=- r^{21}(v,u)$,
\item[{\it (ii)}]
the function $r(u,v)$ depends only on the difference $z:=u-v$.
\end{itemize}
For the expression (\ref{CYBE2}) these properties are equivalent to
\begin{equation}
r(z)=\sum_{i,j}f_{ij}(z)\,e_i\otimes e_j~,
\label{CYBE4}
\end{equation}
where the 
functions $f_{ij}(z)$ satisfy the relation
\begin{equation}
f_{ij}(z)=-f_{ji}(-z)~.
\label{CYBE5}
\end{equation}

In 1982 A. Belavin and V. Drinfeld  \cite{BD} classified such
$z$-dependent solutions of CYBE (\ref{CYBE1}) provided that $r(z)$ is
a "non-degenerate" meromorphic function. They showed that
there exist three type of the solutions: rational, trigonometric and
elliptic. Namely, if $f_{ij}(z)$ are rational functions of the variable
$z$ then the corresponding solution (\ref{CYBE4}) is called rational,
if $f_{ij}(z)=\phi_{ij}(\exp(kz))$, where $\phi_{ij}(z)$ are rational
functions, then the corresponding  solution (\ref{CYBE4}) is called
trigonometric and, finally, if the functions $f_{ij}(z)$ are elliptic
the solution is elliptic. The elliptic solution exists only for the
case $\mathfrak{g}=\mathfrak{sl}_n$.

The classical Yang-Baxter equation (\ref{CYBE1}) has a quantum analog,
the quantum Yang-Baxter equation (QYBE), ($r(z)\rightarrow R(z)$):
\begin{equation}
\begin{array}{rcl}
&&R^{12}(z_1\!-\!z_2)\,R^{13}(z_1\!-\!z_3)\,R^{23}(z_2\!-\!z_3)\;=
\\[7pt]
&&\phantom{aaaaaaaaaaaaa}
=\;R^{23}(z_2\!-\!z_3)\,R^{13}(z_1\!-\!z_3)\,R^{12}(z_1\!-\!z_2)~,
\label{CYBE6}
\end{array}
\end{equation}
\begin{eqnarray}
R(z)\!\!&=&\!\! 1+ h r(z)+ \mathfrak{O}(h^2)~.\;\;\quad{}
\label{CYBE7}
\end{eqnarray}
Where $h$ is a deformation parameter. Every solution $r(z)$ corresponds
to a quantum deformation of $U(\mathfrak{g}[u])$ - a quantum group.
The simplest rational $r$-matrix corresponds to a quantum group named
Yangian (introduced by Drinfeld in 1985, \cite{D1}). The simplest
trigonometric solution of CYBE corresponds to a $q$-deformation of
$U(\mathfrak{g}[u])$, the elliptic $r$-matrix corresponds to an elliptic
quantum group.

\setcounter{equation}{0}
\section{A rational-trigonometric solution of CYBE}
It turns out that a sum of the simplest rational and trigonometric
$r$-matrices is also a solution of CYBE. For the sake of simplicity we
take the case $\mathfrak{g}=\mathfrak{sl}_2:=\mathfrak{sl}(2,\mathbb C)$.
The simplest rational $r_{rt}(u,v)$ and
trigonometric $r_{tr}(u,v)$ solutions of CYBE have the form\footnote{In
contrast to (\ref{CYBE7}), here in (\ref{RT1}) the deformation parameters
$\eta$ and $\hbar$ are introduced explicitly in the $r$-matrices.}
\begin{equation}
\begin{array}{rcl}
r_{rt}(u-v)\!\!&=\!\!&\displaystyle{\frac{\eta}{u-v}}\;\Omega^{}_2,
\\[11pt]
r_{tr}(u-v)\!\!&=\!\!&\hbar\left(\displaystyle{\frac{\exp{(u-v)}+1}
{\exp{(u-v)}-1}}\;
\Omega^{}_2+e_{-\alpha}\otimes e_{\alpha}-
e_{a}\otimes e_{-\alpha}\right).
\end{array}
\label{RT1}
\end{equation}
Here $\{h_{\alpha}, e_{\pm\alpha}\}$ is a standard Cartan-Weyl basis
of $\mathfrak{sl}_2$ with the defining relations:
\begin{equation}
[e_{\alpha},\,e_{-\alpha}]=h_{\alpha}~,\qquad
[h_{\alpha},\,e_{\pm\alpha}]=\pm(\alpha,\alpha)e_{\pm\alpha}
\qquad ((\alpha,\alpha)=2)~.
\label{RT2}
\end{equation}
The element $\Omega^{}_2$ is called the Casimir two-tensor:
\begin{equation}
\Omega^{}_2=\mbox{\Large$\frac{1}{2}$}\Big(\Delta_0(C_2)-
C_2\otimes1-1\otimes C_2\Big)~,
\label{RT3}
\end{equation}
where
\begin{equation}
C_2=\mbox{\Large$\frac{1}{2}$}\,h_{\alpha}^2+
e_{-\alpha}e_{\alpha}+e_{\alpha}e_{-\alpha}
\label{RT4}
\end{equation}
is a second order Casimir element, and the symbol $\Delta_0$ is the primitive
co-multipli\-cation
\begin{equation}
\Delta_0(x)=x\otimes1+1\otimes x\qquad(\forall x\in\mathfrak{sl}_2)~.
\label{RT5}
\end{equation}
In the trigonometric $r$-matrix we make change of the variables
$\exp u\to u$, $\exp v\to v$ then the trigonometric $r$-matrix
takes the form
\begin{equation}
r_{tr}(u,v)=\hbar\left(\frac{u+v}{u-v}\;\Omega^{}_2+
e_{-\alpha}\otimes e_{\alpha}-
e_{a}\otimes e_{-\alpha}\right).
\label{RT6}
\end{equation}
This function already does not depend on the difference $u-v$, but it
also satisfies CYBE. Thus, we have to abandon the condition {\it(ii)}
for $r(u,v)$. However, right now it is easy to see that the sum of the
rational and trigonometric solutions of CYBE is also a solution of CYBE.
Indeed, we have
\begin{equation}
r_{tr}(u+\mbox{\Large$\frac{a}{2}$},v+\mbox{\Large$\frac{a}{2}$})=
r_{tr}(u,v)+\mbox{\Large$\frac{\hbar\,a}{\eta}$}\,r_{rt}(u,v)
\qquad{\rm for}\;\;a\in{\mathbb C}~.
\label{RT7}
\end{equation}
Since $r_{tr}(u,v)$ for any $u$ and $v$ satisfies CYBE
\begin{equation}
[r_{tr}^{12}(u_1,u_2),\,r_{tr}^{13}(u_1,u_3)+r_{tr}^{23}(u_2,u_3)]+
[r_{tr}^{13}(u_1,u_3),\,r_{tr}^{23}(u_2,u_3)]=0
\label{RT8}
\end{equation}
the sum $r_{rtr}(u,v):=r_{rt}(u,v)+r_{tr}(u,v)$ also satisfies CYBE.

Each of the classical $r$-matrices $r_{rt}(u,v)$, $r_{tr}(u,v)$ and
$r_{rtr}(u,v)$ defines its own Lie bialgebra ($\mathfrak{sl}_2[u],\delta_i$)
\cite{D3}, where $\delta_i$ is the cocommutator given by
$\delta_i(x)=[x\otimes 1+1\otimes x,r_i]$, ($x\in\mathfrak{sl}_2[u]$,
$r_i$ is one of our classical r-matrices). Each of these Lie bialgebras
corresponds to a quantum deformation of $U(\mathfrak{g}[u])$.
Using the comultiplication formulas of the quantum algebras
$U_{q}(\mathfrak{sl}_2[u])$ and  $Y_{\eta}(\mathfrak{sl}_2)$ one can show
that the bialgebras for $r_{rt}(u,v)$ and $r_{tr}(u,v)$ coincide with
semiclassical limits of $U_{q}(\mathfrak{sl}_2[u])$ and
$Y_{\eta}(\mathfrak{sl}_2)$, i.e., for example,
$\delta_{tr}=\Delta^{\!op}_{q}-\Delta_{q}$ $\mbox{(\rm mod}\,\hbar^{2})$
where $\Delta_{q}^{\!op}$ is a comultiplication opposite to $\Delta_{q}$,
and analogously for $\delta_{rt}$.

The bialgebra with the rational-trigonometric $r$-matrix
$r_{rtr}(u,v):=r_{rt}(u,v)+r_{tr}(u,v)$ corresponds to a two-parameter
deformation of $U(\mathfrak{g}[u])$, which we call  Drinfeldian.
Thus, the Drinfeldian $D_{q\eta}(\mathfrak{sl}_2)$ is a quantization of
$U(\mathfrak{sl}_2[u])$ in direction of the classical $r$-matrix $r_{rtr}(u,v)$
and moreover the quantum algebra $U_{q}(\mathfrak{sl}_2[u])$ and the Yangian
$Y_{\eta}(\mathfrak{sl}_2)$ are limit quantum algebras of
$D_{q\eta}(\mathfrak{sl}_2)$ when the parameters of deformation $\eta$ goes
to $0$ and $q$ goes to $1$, correspondingly.

Below we present a very simple scheme for construction of defining relations
of the Drinfeldian $D_{q\eta}(\mathfrak{g})$ for any finite-dimensional
complex simple Lie algebra $\mathfrak{g}$ using the defining relations of the
quantum ($q$-deformed) algebra $U_{q}(\mathfrak{g}[u])$, and 
a limited pass $q\to1$ gives us very simple defining relations for the Yangian
$Y_{\eta}(\mathfrak{g})$. Thus, we find very simple connection between the
quantum affine algebra $U_{q}(\mathfrak{g}[u])$ and the Yangian
$Y_{\eta}(\mathfrak{g})$.

In order to compare our realization of the Yangian in terms of a Chevalley
basis with Drinfeld's realizations at first we present the Drinfeld
realizations in an explicit form. There are three realizations  of
$Y_{\eta}(\mathfrak{g})$, given by Drinfeld \cite{D1} - \cite{D3}.
We give two of them.

\setcounter{equation}{0}
\section{Drinfeld's realizations of Yangians}
Let $\mathfrak{g}$ be a finite-dimensional complex simple Lie algebra.
Fix a non-zero invariant bilinear form $(\;,\;)$ on $\mathfrak{g}$, and let
$\{I_\alpha\}$ be an orthonormal basis of $\mathfrak{g}$ with respect to
$(\;,\;)$.\\

{\large\sf 1). The first Drinfeld realization of Yangians}.
\begin{definition}\cite{D1}
The Yangian $Y_{\eta}(\mathfrak{g})$ is generated as an associative algebra
over ${\mathbb C}[[\eta]]$ by the Lie algebra $\mathfrak{g}$ and elements
$J(x)$, $x\in\mathfrak{g}$, with the defining relations:
\begin{equation}
\begin{array}{rcl}
{}\quad J(\lambda x+\mu y)\!\!&=\!\!&\lambda J(x)+\mu J(y),
\\[5pt]
J([x,y])\!\!&=\!\!&[x,J(y)]
\qquad{\rm for}\;\;x,y\in\mathfrak{g},\;\;\lambda,\mu\in{\mathbb C}~,
\end{array}
\label{Y1}
\end{equation}
\begin{equation}
\begin{array}{rcl}
{\rm if}\;\;\;\displaystyle{\sum\limits_i}
[x_i,y_i]\!\!&=\!\!&0\quad{\rm for}\;\,x_i,y_i\in\mathfrak{g}
\quad\Rightarrow\;\;
\\
\displaystyle{\sum\limits_i}[J(x_i),\,J(y_i)]\!\!&=\!\!&
\displaystyle{\frac{\eta^2}{12}
\sum\limits_i\sum\limits_{\alpha,\beta,\gamma}}([[x_i,I_{\alpha}],
[y_i,I_{\beta}]],I_{\gamma})
\{I_{\alpha},I_{\beta},I_{\gamma}\}~,
\label{Y2}
\end{array}
\end{equation}%
\begin{equation}
\begin{array}{rcl}
{\rm if}\;\;\;\displaystyle{\sum\limits_i}
[[x_i,y_i],z_i]\!\!&=\!\!&0 \quad{\rm for}\;\,x_i,y_i,z_i\in\mathfrak{g}
\quad\Rightarrow\;\;
\\
\displaystyle{\sum_i}[[J(x_i),J(y_i)],J(z_i)]\!\!&=\!\!&
\displaystyle{\frac{\eta^2}{4}
\sum\limits_i\!\sum\limits_{\alpha,\beta,\gamma}}\!\!
f(x_i,y_i,z_i,I_{\alpha},I_{\beta},I_{\gamma})
\\
\!\!&\!\!&\phantom{aaaaaaaaaa}\times\,
\{I_{\alpha},I_{\beta},J(I_{\gamma})\}~,
\end{array}
\label{Y3}
\end{equation}
where\footnote{\rm In the case $\mathfrak{g}=\mathfrak{sl}_2(\mathbb C)$
there is a more complicated relation instead (\ref{Y2}) (see \cite{D1}).}
the notations are used $\{a_1,a_2,a_3\}:=(1/6)\!\!\!\sum\limits_{i\neq j\neq k}
\!\!\!\!a_ia_ja_k$ and $f(x,y,z,a,b,c):=\mathop{\mbox{\rm Alt}}\limits_{x,y}\,
\mathop{\mbox{\rm Sym}}\limits_{x,z}\left([x,[y,a]],[[z,b],c]\right)$.
\noindent
A comultiplication map ($\Delta_{\eta}:\,Y_{\eta}(\mathfrak{g})\rightarrow
Y_{\eta}(\mathfrak{g})\otimes Y_{\eta}(\mathfrak{g})$), an antipode
($S_{\eta}:\,Y_{\eta}(\mathfrak{g})\rightarrow Y_{\eta}(\mathfrak{g})$) and
a counite ($\varepsilon_{\eta}:\,Y_{\eta}(\mathfrak{g})\rightarrow {\mathbb C}$)
are given by the formulas ($x\in\mathfrak{g}$)
\begin{equation}
\begin{array}{rcl}
\Delta_{\eta}(x)\!\!&=\!\!&x\otimes1+1\otimes x~,
\\[3pt]
{}\qquad\qquad
\Delta_{\eta}(J(x))\!\!&=\!\!&J(x)\otimes1+1\otimes J(x)+
\displaystyle{\frac{\eta}{2}}\,[x\otimes
1,\Omega_2^{}]~,\;\,\phantom{aaaaaaaaa}
\end{array}
\label{Y4}
\end{equation}
\begin{eqnarray}
S_{\eta}(x)\!\!&=\!\!&-x~,\quad\;
S_{\eta}(J(x))\;=-J(x)+\displaystyle{\frac{\eta}{4}}\,\lambda\,x~,
\label{Y5}
\\[3pt]
{}\quad\;\;\varepsilon_{\eta}(x)\!\!&=\!\!&\varepsilon_{\eta}(J(x))=0~,
\quad\;\;
\varepsilon_{\eta}(1)\;=\;1~,
\label{Y5'}
\end{eqnarray}
where $\Omega_2^{}$ is the Casimir two-tensor
($\Omega_2^{}=\sum\limits_{\alpha}I_{\alpha}\otimes I_{\alpha}$),
$\lambda$ is the eigenvalue of the Casimir operator
$C_2=\sum\limits_{\alpha}I_{\alpha}I_{\alpha}$ in the adjoint
representation of $\mathfrak{g}$ in $\mathfrak{g}$.
\label{YD1}
\end{definition}
We may specialize the formal parameter $\eta$ to any complex number
$\nu\in\mathbb C$, however, the resulting Hopf algebra $Y_{\nu}(\mathfrak{g})$
(over $\mathbb C$) is essentially independent of $\nu$, provided that
$\nu\neq0$. It means that any two Hopf algebras $Y_{\nu}^{}(\mathfrak{g})$ and
$Y_{\nu'}(\mathfrak{g})$ with $\nu\neq\nu';\;\nu,\nu'\neq 0$ are isomorphic.
Thus, we can as well take $\nu=1$ and therefore the parameter $\eta$ is
usually dropped. However, for convenience of  passage to the limit
$\eta\to 0$ we shall remain the formal parameter $\eta$.\\[0pt]

{\large\sf 2). The second Drinfeld realization of Yangians}.
Let $A=(a_{ij})_{i,j=1}^l$ be a standard Cartan matrix of $\mathfrak{g}$,
$\Pi:= \{\alpha_1,\ldots, a_l\}$ be a system of simple roots ($l$ is rank of
$\mathfrak{g}$), and let $B_{ij}:=\frac{1}{2}(\alpha_i,\alpha_j)$.
\begin{theorem}\cite{D2}
The Yangian $Y_{\eta}(\mathfrak{g})$ is isomorphic to the associative
algebra over ${\mathbb C}[[\eta]]$ with the generators:
\begin{equation}
\xi_{in}^{+}~,\quad\xi_{in}^{-}~,\quad\varphi_{in}^{}\qquad {\rm for}
\quad i=1,2,\ldots,l;\;\;n=0,1,2,\ldots,
\label{Y6}
\end{equation}
and the following defining relations:
\begin{eqnarray}
[\varphi_{in}^{},\,\varphi_{jm}^{}]\!\!&=\!\!&0~,
\label{Y7}
\\[3pt]
[\varphi_{i0}^{},\,\xi_{jn}^{\pm}]\!\!&=\!\!&
\pm 2B_{ij}\,\xi_{jn}^{\pm}~,
\label{Y8}
\\[3pt]
[\xi_{in}^{+},\,\xi_{jm}^{-}]\!\!&=\!\!&
\delta_{ij}\,\varphi_{j\,n+m}^{}~,
\label{Y9}
\end{eqnarray}
\begin{eqnarray}
[\varphi_{i\,n+1}^{},\,\xi_{jm}^{\pm}]-[\varphi_{i\,n}^{},\,
\xi_{j\,m+1}^{\pm}]\!\!\!&=\!\!\!&\pm\eta\,B_{ij}\left(\varphi_{i\,n}^{}
\xi_{jm}^{\pm}+\xi_{jm}^{\pm}\varphi_{in}^{}\right),
\label{Y10}
\\[5pt]
[\xi_{i\,n+1}^{\pm}\,\xi_{jm}^{\pm}]-[\xi_{i\,n}^{\pm},\,
\xi_{j\,m+1}^{\pm}]\!\!\!&=\!\!\!&\pm\eta\,B_{ij}
\left(\xi_{i\,n}^{\pm}\xi_{jm}^{\pm}
+\xi_{jm}^{\pm}\xi_{in}^{\pm}\right),
\label{Y11}
\\[5pt]
\mathop{\mbox{\rm Sym}}_{n_1^{},n_2^{},\ldots,n_k^{}}
\!\!\!
[\xi_{i\,n_1^{}}^{\pm}\!,[\xi_{in_2^{}}^{\pm}
[\dots[\xi_{i\,n_k^{}}^{\pm}\!,\xi_{jm}^{\pm}]\ldots]]]\!\!&=\!\!&0
\;\;{\rm for}\;\,i\neq j,\;\, k=1-A_{ij}.
\label{Y11'}
\end{eqnarray}
\label{YD2}
\end{theorem}
Explicit formulas for the action of the comultiplication $\Delta_\eta$
on the generates $\xi_{in}^{\pm}$, $\varphi_{in}^{}$ are rather
cumbersome (see \cite{KT3}) and they are not given here.

For the Yangian $Y_{\eta}(\mathfrak{sl}(n,\mathbb C)))$  Drinfeld has
also given a third realization. It is presented in terms of $RLL$-relations
(see details in \cite{D2, MNO, RTF}).

All these realizations of Yangians are not minimal, i.e. they are not given
in terms of a Chevalley basis. However, the minimal realization exists and
we present it here using connection of Yangians with quantum non-twisted
affine algebras.

\setcounter{equation}{0}
\section{Quantum algebra $U_q(\widetilde{\mathfrak{g}[u]})$}
Again, let $\mathfrak{g}$ be a finite-dimensional complex simple Lie algebra
of rank $l$ with the standard Cartan matrix $A=(a_{ij})_{i,j=1}^l$, with
the system of simple roots $\Pi:= \{\alpha_1,\ldots, a_l\}$, and with the
Chevalley basis $h_{\alpha_i}$, $e_{\pm\alpha_i}'$ $(i=1,2,\ldots, l)$.
Let $\theta$ be the maximal positive root of $\mathfrak{g}$. The corresponding
non-twisted affine algebra $\hat{\mathfrak{g}}$ is generated by $\mathfrak{g}$
and the additional affine elements $e_{\pm(\delta-\theta)}'$ and
$h_{\delta-\theta}$. It is well known \cite{Kac} that
$\hat{\mathfrak{g}}\simeq\widetilde{\mathfrak{g}[u,u^{-1}]}$.

Let $\widetilde{\mathfrak{g}[u]}$ be a Lie algebra generated by $\mathfrak{g}$,
the positive root vector $e_{\delta-\theta}'=ue_{-\theta}'$ and the Cartan element
$h_{\delta-\theta}$, i.e. $\widetilde{\mathfrak{g}[u]}\simeq\mathfrak{g}[u]\oplus
{\mathbb C}h_{\delta}$, where $h_{\delta}$ is a central element. The standard
defining relations of $\widetilde{\mathfrak{g}[u]}$ and its universal
enveloping algebra $U(\widetilde{\mathfrak{g}[u]})$ are given by the formulas:
\begin{eqnarray}
[h_{\delta}^{},\mbox{\rm everything}]\!\!&=\!\!&0~,
\label{tolsA1}
\\[3pt]
[h_{\alpha_{i}}^{},h_{\alpha_{j}}^{}]\!\!&=\!\!&0~,
\label{tolsA1'}
\\[3pt]
[h_{\alpha_{i}}^{},e_{\pm\alpha_{j}}']\!\!&=\!\!&
\pm(\alpha_{i},\alpha_{j})e_{\pm\alpha_{j}}'~,
\label{tolsA2}
\\[3pt]
[e_{\alpha_{i}}',e_{-\alpha_{j}}']\!\!&=\!\!&\delta_{ij}h_{\alpha_{i}}~,
\label{tolsA2'}
\\[3pt]
(\ad e_{\pm\alpha_{i}}'\!)^{n_{ij}} e_{\pm\alpha_{j}}'\!\!&=\!\!&0
\quad\;{\rm for}\;\;i\neq j,\;\;n_{ij}\!:=\!1\!-\!a_{ij}~,
\label{tolsA3}
\\[3pt]
[h_{\alpha_i}^{},e_{\delta-\theta}']\!\!&=\!\!&-
(\alpha_i,\theta)\,e_{\delta-\theta}'~,
\label{tolsA3'}
\\[3pt]
[e_{-\alpha_i}',e_{\delta-\theta}']\!\!&=\!\!&0~,
\label{tolsA4}
\\[3pt]
(\ad e_{\alpha_i}')^{n_{i0}}e_{\delta-\theta}'\!\!&=\!\!&0
\quad\;{\rm for}\;\,n_{i0}=1\!+\!2(\alpha_i,\theta)/(\alpha_i,\alpha_i),
\label{tolsA5}
\\[3pt]
[[e_{\alpha_i}',e_{\delta-\theta}'],e_{\delta-\theta}']\!\!&=\!\!&0\quad\;
{\rm for}\;\,\mathfrak{g}\ne\mathfrak{sl}_2 \;\;{\rm and}\;\;(\alpha_i,\theta)\ne 0~,
\label{tolsA6}
\\[3pt]
[[[e_{\alpha}',e_{\delta-\alpha}'],e_{\delta-\alpha}'],e_{\delta-\alpha}']
\!\!&=\!\!&0\quad\;{\rm for}\;\,\mathfrak{g}=\mathfrak{sl}_2~.
\label{tolsA7}
\end{eqnarray}
Here $"\ad"$ is the adjoint action of $\widetilde{\mathfrak{g}[u]}$ in
$\widetilde{\mathfrak{g}[u]}$, i.e. $(\ad x)y=[x,y]$ for
$x,y\in\widetilde{\mathfrak{g}[u]}$. The relations (\ref{tolsA6}) relate to
the case $\mathfrak{g}\ne\mathfrak{sl}_2$, and the relation (\ref{tolsA7})
belongs to the case $\mathfrak{g}=\mathfrak{sl}_2$ (in this case $\theta=\alpha$).

The defining relations of the quantum algebra $U_q(\widetilde{\mathfrak{g}[u]})$,
which is a $q$-defor\-mation of $U(\widetilde{\mathfrak{g}[u]})$, can be obtained
from the defining relations of the quantum algebra $U_q(\hat{\mathfrak{g}})$ by
removing the relations with the negative affine root vector $e_{-\delta+\theta}$.
If $k_{\delta}^{\pm 1}:=q^{\pm h_\delta}$,
$k_{\alpha_i}^{\pm 1}:=q^{\pm h_{\alpha_i}}$, $e_{\pm\alpha_i}$
$(i=1,2,\ldots, l)$, and $e_{\pm(\delta-\theta)}$ is the Chevalley
basis of $U_q(\hat{\mathfrak{g}})$ then the standard defining relations of
$U_q(\widetilde{\mathfrak{g}[u]})$ are given by the formulas:
\begin{eqnarray}
[k_{\delta}^{\pm 1},{\rm everything}]\!\!&=\!\!&0~,
\label{tolsA8}
\\[5pt]
k_{\alpha_i}^{\pm 1}k_{\alpha_j}^{\pm 1}\!\!&=\!\!&
k_{\alpha_j}^{\pm 1}k_{\alpha_i}^{\pm 1}~,
\label{tolsA8'}
\\[5pt]
k_{\delta}^{}k^{-1}_{\delta}=k^{-1}_{\delta}k_{\delta}^{}\!\!&=\!\!&
k_{\alpha_i}^{}k^{-1}_{\alpha_i}=k^{-1}_{\alpha_i}k_{\alpha_i}^{}=1~,
\label{tolsA9}
\\[5pt]
k_{\alpha_i}^{}e_{\pm\alpha_j}^{}k^{-1}_{\alpha_i}\!\!&=\!\!&
q^{\pm(\alpha_i,\alpha_j)}e_{\pm\alpha_j}^{}~,
\label{tolsA9'}
\\[5pt]
[e_{\alpha_i}^{},e_{-\alpha_i}^{}]\!\!&=\!\!&
\frac{k_{\alpha_i}-k_{\alpha_i}^{-1}}{q-q^{-1}}~,
\label{tolsA10}
\\[5pt]
(\ad_{q}e_{\pm\alpha_{i}}^{}\!)^{n_{ij}}e_{\pm\alpha_{j}^{}}\!\!&=\!\!&0
\quad{\rm for}\;\,i\neq j,\;\;n_{ij}\!:=\!1\!-\!a_{ij}~,
\label{tolsA11}
\\[5pt]
k_{\alpha_i}^{}e_{\delta-\theta}^{}k^{-1}_{\alpha_i}\!\!&=\!\!&
q^{-(\alpha_i,\theta)}e_{\delta-\theta}^{}~,
\label{tolsA11'}
\\[5pt]
[e_{-\alpha_i}^{},e_{\delta-\theta}^{}]\!\!&=\!\!&0~,
\label{tolsA12}
\\[5pt]
(\ad_{q}e_{\alpha_i}^{})^{n_{i0}}e_{\delta-\theta}^{}\!\!&=\!\!&0
\quad{\rm for}\;\;n_{i0}=1\!+\!2(\alpha_i,\theta)/(\alpha_i,\alpha_i),
\label{tolsA13}
\\[5pt]
[[e_{\alpha_i}^{},e_{\delta-\theta}^{}]_{q},
e_{\delta-\theta}^{}]_{q}\!\!&=\!\!&0\quad{\rm for}\;\;
\mathfrak{g}\ne\mathfrak{sl}_2\;\;{\rm and}\;\;(\alpha_i,\theta)\ne0~,
\label{tolsA14}
\\[5pt]
[[[e_{\alpha}^{},e_{\delta-\alpha}^{}]_{q},
e_{\delta-\alpha}^{}]_{q},e_{\delta-\alpha}^{}]_{q}\!\!&=\!\!&0
\quad{\rm for}\;\;\mathfrak{g}=\mathfrak{sl}_2~,
\label{tolsA15}
\end{eqnarray}
where $(\ad_{q}e_{\beta})e_{\gamma}$ is the $q$-commutator:
\begin{equation}
\quad\quad\;{}
(\ad_qe_{\beta}^{})e_{\gamma}^{}:= [e_\beta^{},e_\gamma^{}]_q:=
e_\beta^{}e_\gamma^{}-q^{(\beta,\gamma)}e_\gamma^{}e_\beta^{}~.
\label{tolsA16}
\end{equation}
The comultiplication $\Delta_{q}$, the antipode $S_{q}$, and
the counite $\varepsilon_{q}$ of $U_{q}(\widetilde{\mathfrak{g}[u]})$
are given by
\begin{equation}
\begin{array}{rcl}
\Delta_{q}(k_{\gamma}^{\pm 1})\!\!&=\!\!&
k_{\gamma}^{\pm 1}\otimes k_{\gamma}^{\pm 1}
\quad\;(\gamma=\delta,\;\alpha_i)~,
\\[5pt]
\Delta_{q}(e_{-\alpha_i}^{})\!\!&=\!\!&
e_{-\alpha_i}^{}\otimes k_{\alpha_i}^{}+1\otimes
e_{-\alpha_i}^{}~,
\\[5pt]
\Delta_{q}(e_{\alpha_i}^{})\!\!&=\!\!&
e_{\alpha_i}^{}\otimes1+k_{\alpha_i}^{-1}\otimes
e_{\alpha_i}^{}~,\qquad
\\[5pt]
\Delta_{q}(e_{\delta-\theta}^{})\!\!&=\!\!&e_{\delta-\theta}^{}
\otimes1+k_{\delta-\theta}^{-1}\otimes
e_{\delta-\theta}^{}~,
\label{tolsA19}
\end{array}
\end{equation}
\begin{equation}
\begin{array}{rcl}
S_{q}(k_{\gamma}^{\pm 1})\!\!&=\!\!&k_{\gamma}^{\mp 1}
\qquad(\gamma=\delta,\;\alpha_i)~,\phantom{aaaaa}
\\[5pt]
S_{q}(e_{-\alpha_i}^{})\!\!&=\!\!&-e_{-\alpha_i}^{}k_{\alpha_i}^{-1}~,
\\[5pt]
S_{q}(e_{\alpha_i}^{})\!\!&=\!\!&-k_{\alpha_i}^{}e_{\alpha_i}^{}~,
\\[5pt]
S_{q}(e_{\delta-\theta}^{})\!\!&=\!\!&
-k_{\delta-\theta}^{}e_{\delta-\theta}^{}~,
\end{array}
\end{equation}
\begin{equation}
\varepsilon_{q}(e_{\pm\alpha_i})~=~\varepsilon_{q}(e_{\delta-\theta})=0~,
\qquad\varepsilon_{q}(k_{\delta}^{\pm 1})
~=~\varepsilon_{q}(k_{\alpha_i}^{\pm 1})~=~1~.
\label{tolsA22'}
\end{equation}
We set
\begin{equation}
\qquad\;{}k_{\delta-\theta}^{}=\,k_{\delta}^{}k_{\alpha_1}^{-n_1}
k_{\alpha_2}^{-n_2}\cdots k_{\alpha_l}^{-n_l}
\label{tolsA23}
\end{equation}
if $\theta=n_1^{}\alpha_1^{}+n_2^{}\alpha_2^{}+\cdots\,+ n_l^{}\alpha_l^{}$.
\par
It is easy to check that the following proposition is valid.
\begin{proposition}
There is a one-parameter group of Hopf algebra automorphisms
${\mathfrak T}_a$ of $U_q(\widetilde{\mathfrak{g}[u]})$, $a\in {\mathbb C}$,
given by
\begin{equation}
{\mathfrak T}_a(k_{\gamma}^{})\,=\,k_{\gamma}^{}~,\qquad
{\mathfrak T}_a(e_{\pm\alpha_i})\,=\,e_{\pm\alpha_i}~,\qquad
{\mathfrak T}_a(e_{\delta-\theta}^{})\,=\,a\,e_{\delta-\theta}^{}~.
\quad{}
\label{tolsA24}
\end{equation}
\label{tolsPA1}
\end{proposition}
Now starting from the quantum algebra $U_q(\widetilde{\mathfrak{g}[u]})$
we construct the two-param\-eter Hopf algebra which coincides with the Yangian
$Y_{\eta}(\mathfrak{g})$ at $q=1$.

\setcounter{equation}{0}
\section{Drinfeldian $D_{q\eta}(\mathfrak{g})$}

The quantum algebra $U_{q}(\widetilde{\mathfrak{g}[u]})$ does not contain
any singular elements when $q$ goes to 1 $(q=\exp(\hbar))$. This means
that $\lim\limits_{q\to1}x\in U(\widetilde{\mathfrak{g}[u]})$ for any
$x\in U_q(\widetilde{\mathfrak{g}[u]})$. Now we extend
$U_q(\widetilde{\mathfrak{g}[u]})$ by singular elements,
i.e. we consider the algebra
\begin{equation}
\bar{U}_{q}(\widetilde{\mathfrak{g}[u]}):=U_q(\widetilde{\mathfrak{g}[u]})
\otimes_{{\mathbb C}[[\hbar]]}{\mathbb C}((\hbar))~.
\label{tolsD1}
\end{equation}
This algebra contains singular elements, i.e. elements which have not
any limit when $q\to1$. For example, if
$x\in U_q(\widetilde{\mathfrak{g}[u]})$ and $\lim\limits_{q\to1}x\ne0$
then the element $x/(q-q^{-1})$ is singular.

Let $\tilde{e}_{-\theta}^{}\in U_q(\widetilde{\mathfrak{g}[u]})$ be any
element of the weight $-\theta$, i.e.
\begin{equation}
k_{\alpha_i}^{}\tilde{e}_{-\theta}^{}\,k^{-1}_{\alpha_i}=
q^{-(\alpha_i,\theta)}\tilde{e}_{-\theta}^{}~,\quad{}
\label{tolsD2}
\end{equation}
and, moreover,
\begin{equation}
\lim\limits_{q\to 1}\tilde{e}_{-\theta}^{}=e_{-\theta}'~,\qquad\quad
\label{tolsD3}
\end{equation}
where $e_{-\theta}'$ is a root vector of $\mathfrak{g}$ with the minimal
weight $-\theta$. Let us introduce the new affine generator
\begin{equation}
\quad{}\xi_{\delta-\theta}^{}=e_{\delta-\theta}^{}+
\frac{\eta\phantom{a}}{q-q^{-1}}
\,\tilde{e}_{-\theta}^{}~.
\label{tolsD4}
\end{equation}
Then the defining relations (\ref{tolsA11'})--(\ref{tolsA15})
take the form
\begin{eqnarray}
k_{\alpha_i}^{}\xi_{\delta-\theta}^{}k^{-1}_{\alpha_i}\!\!&=\!\!&
q^{-(\alpha_i,\theta)}\xi_{\delta-\theta}^{}~,
\label{tolsD5}
\\[3pt]
[e_{-\alpha_i}^{},\xi_{\delta-\theta}^{}]\!\!&=\!\!&
\tau\,[e_{-\alpha_i}^{},\tilde{e}_{-\theta}^{}]~,
\label{tolsD6}
\\[3pt]
({\rm ad}_{q}e_{\alpha_i}^{})^{n_{i0}^{}}\xi_{\delta-\theta}^{}\!\!&=\!\!&
\tau\,({\rm ad}_{q}e_{\alpha_i}^{})^{n_{i0}^{}}\tilde{e}_{-\theta}^{}
\label{tolsD7}
\end{eqnarray}
for $n_{i0}^{}=1+2(\alpha_i,\theta)/(\alpha_i,\alpha_i)$~, and
\begin{equation}
\begin{array}{rcl}
[[e_{\alpha_i}^{},\xi_{\delta-\theta}^{}]_{q},
\xi_{\delta-\theta}^{}]_{q}\!\!&=\!\!&-\tau^{2}[[e_{\alpha_i}^{},
\tilde{e}_{-\theta}^{}]_{q},\tilde{e}_{-\theta}^{}]_{q}
\\[3pt]
&&+\tau\Big([[e_{\alpha_i}^{},\tilde{e}_{-\theta}^{}]_{q},
\xi_{\delta-\theta}]_{q}+
[[e_{\alpha_i}^{},\xi_{\delta-\theta}^{}]_{q},
\tilde{e}_{-\theta}^{}]_{q}\Big)
\end{array}
\label{tolsD8}
\end{equation}
for $\mathfrak{g}\ne\mathfrak{sl}_2^{}$ and $(\alpha_i,\theta)\ne 0$,
\begin{equation}
\begin{array}{rcl}
&&[[[e_{\alpha}^{},\xi_{\delta-\alpha}^{}]_{q},
\xi_{\delta-\alpha}^{}]_{q},\xi_{\delta-\alpha}^{}]_{q}
=\tau\Big([[[e_{\alpha}^{},\tilde{e}_{-\alpha}^{}]_{q},
\xi_{\delta-\alpha}^{}]_{q},\xi_{\delta-\alpha}^{}]_{q}
\\[3pt]
&&\qquad\qquad{}
+[[[e_{\alpha}^{},\xi_{\delta-\alpha}^{}]_{q},
\tilde{e}_{-\alpha}^{}]_{q},\xi_{\delta-\alpha}^{}]_{q}
+[[[e_{\alpha}^{},\xi_{\delta-\alpha}^{}]_{q},
\xi_{\delta-\alpha}^{}]_{q},\tilde{e}_{-\alpha}^{}]_{q}\Big)
\\[3pt]
&&\qquad\qquad{}
-\tau^{2}\Bigl([[[e_{\alpha}^{},\tilde{e}_{-\alpha}^{}]_{q},
\tilde{e}_{-\alpha}^{}]_{q},\xi_{\delta-\alpha}^{}]_{q}
+[[[e_{\alpha}^{},\tilde{e}_{-\alpha}^{}]_{q},
\xi_{\delta-\alpha}^{}]_{q},\tilde{e}_{-\alpha}^{}]_{q}
\\[3pt]
&&\qquad\qquad{}
+[[[e_{\alpha}^{},\xi_{\delta-\alpha}^{}]_{q},
\tilde{e}_{-\alpha}^{}]_{q},\tilde{e}_{-\alpha}^{}]_{q}\Big)
+\tau^{3}[[[e_{\alpha}^{},\tilde{e}_{-\alpha}^{}]_{q},
\tilde{e}_{-\alpha}^{}]_{q},\tilde{e}_{-\alpha}^{}]_{q}
\label{tolsD9}
\end{array}
\end{equation}
for $\mathfrak{g}=\mathfrak{sl}_2^{}$. Here and elsewhere
$\tau:=\eta/(q-q^{-1})$.

The formulas (\ref{tolsA19}) for the comultiplication $\Delta_{q\eta}$
and the antipode $S_{q\eta}$ are rewritten as follows
\begin{equation}
\Delta_{q\eta}(x)=\Delta_{q}(x)
\quad\;(x\in U_{q}(g)\otimes\mathbb{C}k_{\delta}^{}),\qquad{}
\label{tolsD10}
\end{equation}
\begin{equation}
\begin{array}{rcl}
\Delta_{q\eta}(\xi_{\delta-\theta}^{})\!\!&=\!\!&
\xi_{\delta-\theta}^{}\otimes 1+k_{\delta-\theta}^{-1}
\otimes\xi_{\delta-\theta}^{}
\\[5pt]
&&{}+\tau\Bigl(\Delta_{q}(\tilde{e}_{-\theta}^{})
-\tilde{e}_{-\theta}^{}\otimes 1-
k_{\delta-\theta}^{-1}\otimes\tilde{e}_{-\theta}^{}\Bigr),
\end{array}
\label{tolsD11}
\end{equation}
\begin{eqnarray}
S_{q\eta}(x)\!\!&=\!\!&S_{q}(x)
\quad\;(x\in U_{q}(\mathfrak{g})\otimes\mathbb{C}k_{\delta}^{}),
\label{tolsD12}
\\[5pt]
S_{q\eta}(\xi_{\delta-\theta}^{})\!\!&=\!\!&
-k_{\delta-\theta}^{}\xi_{\delta-\theta}^{}
+\tau\Bigl(S_{q}(\tilde{e}_{-\theta}^{})+
k_{\delta-\theta}^{}\tilde{e}_{-\theta}^{}\Bigr).
\label{tolsD13}
\end{eqnarray}
The following theorem is valid.
\begin{theorem}
If the element $\tilde{e}_{-\theta}\in U_q(\widetilde{\mathfrak{g}[u]})$
satisfies the conditions (\ref{tolsD2}) and (\ref{tolsD3}) then the right-hand
sides of the relations (\ref{tolsD6})--(\ref{tolsD9}) and (\ref{tolsD11}),
(\ref{tolsD13}) are nonsingular at $q=1$.
\label{tolsDT1}
\end{theorem}
Thus, the associative algebra $D_{q\eta}(\mathfrak{g})$
over $\mathbb C[[\hbar,\eta]]$ generated
by the quantum algebra $U_{q}(\mathfrak{g})$,
the central element $k_{\delta}^{\pm 1}$ and by the affine element
$\xi_{\delta-\theta}^{}$ with the relations (\ref{tolsD5})--(\ref{tolsD9}),
where $\tilde{e}_{-\theta}\in U_q(\widetilde{\mathfrak{g}[u]})$ satisfies
the conditions (\ref{tolsD2}) and (\ref{tolsD3}), and with the Hopf structure
given by the formulas (\ref{tolsD10})--(\ref{tolsD13}) is nonsingular at
$q=1$.
\begin{theorem}
(i) The Hopf algebra $D_{q\eta}(\mathfrak{g})$ is a two-parameter quantization
of $U(\widetilde{\mathfrak{g}[u]})$, ($\widetilde{\mathfrak{g}[u]}:=
\mathfrak{g}[u]\oplus h_{\delta}$), in the direction of a classical $r$-matrix
which is a sum of the simplest rational and trigonometric $r$-matrices.

\noindent
(ii) The Hopf algebra $D_{q=1,\eta}(\mathfrak{g}):=D_{q\eta}(\mathfrak{g})
/\hbar D_{q\eta}(\mathfrak{g})$ is isomorphic to the Yangian
$Y_{\eta}(\mathfrak{g})$ (with the additional central element $h_{\delta}$).
Moreover, $D_{q\eta=0}(\mathfrak{g})=U_q(\widetilde{\mathfrak{g}[u]})$.
\label{DT2}
\end{theorem}
By virtue of this theorem, the Hopf algebra $D_{q\eta}(\mathfrak{g})$ can be
called as {\it the rational-trigono\-metric quantum algebra} or
{\it the q-deformed Yangian}. 
We also call this algebra as {\it the Drinfeldian} (in honor
of V.G. Drinfeld, who first observed existence of such subalgebra
in the extension (\ref{tolsD1}) of the $q$-deformed affine algebra
$U_q(\hat{\mathfrak{g}})$ \cite{D3}).

\begin{remark}
Since the defining relations for $D_{q\eta}(\mathfrak{g})$
and $U(\widetilde{\mathfrak{g}[u]})$ in terms of the Chevalley basis differ
only in the right-hand sides of the relations (\ref{tolsD6})--(\ref{tolsD9})
and (\ref{tolsA12})--(\ref{tolsA15}), therefore the Dynkin diagram of
$\mathfrak{g}[u]$ can be also used for classification of the Drinfeldian
$D_{q\eta}(\mathfrak{g})$ and the Yangian $Y_{\eta}(\mathfrak{g})$.
\end{remark}
As an immediate consequence of Proposition \ref{tolsPA1} we have
the following result.
\begin{proposition}
There is a one-parameter group of Hopf algebra automorphisms
${\mathfrak T}_a$ of $D_{q\eta}(\mathfrak{g})$, $a\in{\mathbb C}$, given by
\begin{equation}
\begin{array}{rcl}
{\mathfrak T}_a(x)\!\!&=\!\!&x\quad\;(x\in U_{q}(g)\otimes k_{\delta}^{})~,
\\[5pt]
{\mathfrak T}_a(\xi_{\delta-\theta}^{})\!\!&=\!\!&\Bigl(1-(q-q^{-1})\,a\Bigr)
\xi_{\delta-\theta}^{}+\eta\,a\,\tilde{e}_{-\theta}^{}~.
\label{tolsD14}
\end{array}
\end{equation}
\label{tostPD1}
\end{proposition}
The relations between the Drinfeldian $D_{q\eta}(\mathfrak{g})$ and
the algebras $U_{q}(\mathfrak{g}[u])$, $Y_{\eta}(\mathfrak{g})$,
$U(\mathfrak{g}[u])$ (and also their subalgebras) are shown in the picture:
\\
\begin{eqnarray}
\mbox{\begin{picture}(50,40)
\put(-28,25){$D_{q\eta}(\mathfrak{g})$}
\put(-26,-35){\mbox{$Y_{\eta}(\mathfrak{g})$}}
\put(78,25){\mbox{$U_{q}(\mathfrak{g}[u])$}}
\put(78,-35){\mbox{$U(\mathfrak{g}[u])$}}
\put(23,27){\vector(1,0){40}}
\put(20,-33){\vector(1,0){42}}
\put(-22,10){\vector(0,-1){25}}
\put(85,10){\vector(0,-1){25}}
\put(-68,25){\mbox{$U_{q}(\mathfrak{g})\,\subset$}}
\put(-65,-35){\mbox{$U(\mathfrak{g})\,\subset$}}
\put(117,25){\mbox{$\supset\;U_{q}(\mathfrak{g})$}}
\put(114,-35){\mbox{$\supset\;U(\mathfrak{g})$}\,\,.}
\put(-58,-5){\mbox{\sms$q\!\to\!1$}}
\put(95,-5){\mbox{\sms$q\!\to\!1$}}
\put(28,34){\mbox{\sms$\eta\!\to\!0$}}
\put(28,-26){\mbox{\sms$\eta\!\to\!0$}}
\end{picture}}
\label{tolsDY46}
\nonumber
\end{eqnarray}
\nopagebreak
\\[15pt]
\begin{quote}
\begin{center}
{\fs Fig. 1. Diagram of the limit Hopf algebras
of the Drinfeldian $D_{q\eta}(\mathfrak{g})$\\
and their subalgebras. The arrows show passages to the limits.}
\end{center}
\end{quote}

\setcounter{equation}{0}
\section{Drinfeldians and Yangians for Lie algebras $A_{l-1}^{}$
($l\geq3$)}

An explicit description of the Drinfeldian  $D_{\eta\,q}(\mathfrak{g})$
and Yangian $Y_{\eta}(\mathfrak{g})$ for the simple Lie algebra
$\mathfrak{g}=A_1^{}\simeq\mathfrak{sl}_2:=\mathfrak{sl}(2,\mathbb C)$
was given in \cite{T2, T3} and, therefore, we consider here the case
$\mathfrak{g}=A_{l-1}^{}$, where $l\geq3$.

We first recall the defining relations of the $q$-deformed universal
enveloping algebra $U_q(\mathfrak{sl}_{l})$ ($\mathfrak{sl}_{l}^{}:=
\mathfrak{sl}(l,\mathbb C)\simeq A_{l-1}^{}$) and structure of its
Cartan-Weyl basis.

Let $\Pi:=\{\alpha_{1}^{},\ldots,\alpha_{l-1}^{}\}$ be a system of simple
roots of $\mathfrak{sl}_{l}^{}$ endowed with the following scalar product:
$(\alpha_{i}^{},\alpha_{j}^{})=(\alpha_{j}^{},\alpha_{i}^{})$,
$(\alpha_{i}^{},\alpha_{i}^{})=2$, $(\alpha_{i}^{},\alpha_{i+1}^{})=-1$,
$(\alpha_i^{},\alpha_j^{})=0$ ($|i-j|>1$).
The corresponding Dynkin diagram is presented on the picture:
\\
\begin{eqnarray}
\mbox{\begin{picture}(100,10)
\put(-40,0){\circle{5}}
\put(-37,0){\line(1,0){44}}
\put(-43,-15){\fs$\alpha_{1}^{}$}
\put(10,0){\circle{5}}
\put(13,0){\line(1,0){35}}
\put(7,-15){\fs$\alpha_{2}^{}$}
\put(55,0){$\ldots$}
\put(107,0){\line(-1,0){35}}
\put(110,0){\circle{5}}
\put(113,0){\line(1,0){44}}
\put(100,-15){\fs$\alpha_{l-2}^{}$}
\put(160,0){\circle{5}}
\put(157,-15){\fs$\alpha_{l-1}^{}$}
\end{picture}}
\label{A1'}
\nonumber
\end{eqnarray}
\nopagebreak
\\[0pt]
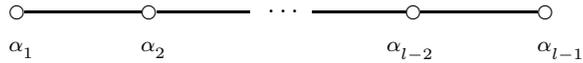
\centerline{\fs Fig. 2. The Dynkin diagram of the Lie algebra
$\mathfrak{sl}_{l}^{}$.}
\\

\noindent
Usually, it is convenient to realize a positive root system $\Delta_{+}^{}$
(with respect to $\Pi$) and a total root system
$\Delta=\Delta_{+}^{}\bigcup(-\Delta_{+}^{})$ of the classical Lie algebras
$A_{l-1}^{}$, $B_l^{}$, $C_l^{}$ and $D_l^{}$ in terms of an orthonormalized
basis of a $l$-dimensional Euclidian space.
Namely, let $\epsilon_i^{}$ ($i=1,2,\ldots, l$) be an orthonormalized
basis of a $l$-dimensional Euclidian space $\mathbb R^{l}$:
$(\epsilon_i^{},\epsilon_j^{})=\delta_{ij}^{}$. In the terms of $\epsilon_i^{}$
the systems $\Pi$, $\Delta_{+}^{}$ and $\Delta$ of $\mathfrak{sl}_{l}^{}$ are
presented as follows:
\begin{eqnarray}
\Pi\!\!&=\!\!&\{\alpha_i=\epsilon_i^{}-\epsilon_{i+1}^{}\,|\,
i=1,2,\ldots,l-1\},
\label{A1}
\\[5pt]
\Delta_{+}^{}\!\!\!\!&=\!\!&
\{\epsilon_i^{}-\epsilon_j^{}\,|\,1\le i<j\le l\},
\label{A2}
\\[5pt]
\Delta_{}\!\!&=\!\!&\Delta_{+}\cup\,(-\Delta_{+})\;=\;
\{\epsilon_i^{}-\epsilon_j^{}\,|\,i\neq j;\,i,j=1,2,\ldots,l\}.
\label{A3}
\end{eqnarray}
The root $\theta:=\alpha_1^{}+\alpha_2^{}+\ldots+\alpha_{l}^{}=
\epsilon_{1}^{}-\epsilon_{l}^{}$ is maximal.
We shall use instead of the simple Lie algebra $\mathfrak{sl}_{l}^{}$
its central extension $\mathfrak{gl}_{l}^{}$.

The quantum algebra $U_{q}(\mathfrak{gl}_{l}^{})$ is generated by the
Chevalley elements $e_{i,-i-1}^{}:=e_{\epsilon_i-\epsilon_{i+1}}^{}$,
$e_{i+1,-i}^{}:=e_{\epsilon_{i+1}-\epsilon_{i}}^{}$ $(i=1,2,\ldots,l-1)$,
and $q^{\pm e_{i,-i}}$ $(i=1,2,\ldots,l)$ with the defining
relations\footnote{If $q^{\sum_{i=1}^{l}e_{i,-i}^{} }=1$ then we obtain
the quantum algebra $U_{q}(\mathfrak{sl}_{l}^{}).$}:
\begin{equation}
\begin{array}{rcl}
q^{e_{i,-i}}q^{-e_{i,-i}}\!\!&=\!\!&q^{-e_{i,-i}}q^{e_{i,-i}}=1~,
\\[7pt]
q^{e_{i,-i}}q^{e_{j,-j}}\!\!&=\!\!&q^{e_{j,-j}}q^{e_{i,-i}}~,
\\[5pt]
q^{e_{i,-i}}e_{j,-k}^{}q^{-e_{i,-i}}\!\!&=\!\!&q^{\delta_{i,j}-\delta_{i,k}}
e_{j,-k}^{}\quad(|j-k|=1)~,
\\[5pt]
[e_{i,-i-1}^{},\,e_{j+1,-j}^{}]\!\!&=\!\!&\delta_{ij}\,
\mbox{\ls$\frac{q^{e_{i,-i}-e_{i+1,-i-1}}\,-\,q^{e_{i+1,-i-1}-e_{i,-i}}}
{q\,-\,q^{-1}}$}~,
\\[7pt]
[e_{i,-i-1}^{},\,e_{j,-j-1}^{}]\!\!&=\!\!&0\quad{\rm for}\;\;|i-j|\geq 2~,
\\[7pt]
[e_{i+1,-i}^{},\,e_{j+1,-j}^{}]\!\!&=\!\!&0\quad{\rm for}\;\;|i-j|\geq 2~,
\\[7pt]
[[e_{i,-i-1}^{},\,e_{j,-j-1}^{}]_{q}^{},\,e_{j,-j-1}^{}]_{q}^{}\!\!&=\!\!&0
\quad{\rm for}\;\;|i-j|=1~,
\\[7pt]
[[e_{i+1,-i}^{},\,e_{j+1,-j}^{}]_{q}^{},\,e_{j+1,-j}^{}]_{q}^{}\!\!&=\!\!&0
\quad{\rm for}\;\;|i-j|=1~.
\label{A4}
\end{array}
\end{equation}
The Hopf structure on $U_{q}(\mathfrak{gl}_{l}^{})$ is given by the following
formulas for a comultiplication $\Delta_{q}$, an antipode $S_{q}$, and a counit
$\varepsilon_{q}$:
\begin{equation}
\begin{array}{rcl}
\Delta_{q}(q^{\pm e_{i,-i}})\!\!&=\!\!&
q^{\pm e_{i,-i}}\otimes q^{\pm e_{i,-i}} ~,
\\[5pt]
\Delta_{q}(e_{i,-i-1}^{})\!\!&=\!\!&
e_{i,-i-1}^{}\otimes 1+q^{e_{i+1,-i-1}-e_{i,-i}}\otimes e_{i,-i-1}^{}~,
\\[5pt]
\Delta_{q}(e_{i+1,-i}^{})\!\!&=\!\!&
e_{i+1,-i}^{}\otimes q^{e_{i,-i}-e_{i+1,-i-1}}+1\otimes e_{i+1,-i}^{}~,
\end{array}
\label{A5}
\end{equation}\
\begin{equation}
\begin{array}{rcl}
S_{q}(q^{\pm e_{i,-i}})\!\!&=\!\!&q^{\mp e_{i,-i}}~,
\phantom{aaaaaaaaaaaaaaaaaaaaaaaa}
\\[5pt]
S_{q}(e_{i,-i-1}^{})\!\!&=\!\!&-q^{e_{i,-i}-e_{i+1,-i-1}}\,e_{i,-i-1}^{}~,
\\[5pt]
S_{q}(e_{i+1,-i}^{})\!\!&=\!\!&-e_{i+1,-i}^{}\,q^{e_{i+1,-i-1}-e_{i,-i}}~,
\end{array}
\label{A6}
\end{equation}
\begin{equation}
\begin{array}{rcccl}
\varepsilon_{q}(q^{\pm e_{i,-i}})^{}\!\!&=\!\!&1,\quad\;\;
\varepsilon_{q}(e_{i,-j}^{})
\!\!&=\!\!&0\quad(|i-j|=1).\phantom{aa}
\end{array}
\label{A7}
\end{equation}

For construction of the composite root vectors $e_{i,-j}^{}$ for
$\mid i-j\mid\geq2$ we fix the following normal ordering of the positive
root system $\Delta_{+}^{}$ (see \cite{T1,KT1,KT2})
\begin{equation}
\begin{array}{rcl}
&&(\epsilon_1^{}\!-\epsilon_2^{},\epsilon_1^{}\!-\epsilon_3^{},\ldots,
\epsilon_1^{}\!-\epsilon_l^{}),(\epsilon_{2}^{}\!-\epsilon_{3}^{},
\epsilon_{2}^{}-\epsilon_{4}^{},\ldots,\epsilon_2^{}\!-\epsilon_{l}^{}),
\ldots,
\\[5pt]
&&(\epsilon_{l-3}^{}\!-\epsilon_{l-2}^{},\epsilon_{l-3}^{}-\epsilon_{l-1}^{},
\epsilon_{l-3}^{}\!-\epsilon_{l}^{}),(\epsilon_{l-2}^{}\!-\epsilon_{l-1}^{},
\epsilon_{l-2}^{}\!-\epsilon_{l}^{}),\epsilon_{l-1}^{}\!-\epsilon_{l}^{}~.
\end{array}
\label{A8}
\end{equation}
According to this ordering we set
\begin{equation}
e_{i,-j}^{}:=[e_{i,-k}^{},\,e_{k,-j}^{}]_{q^{-1}},\qquad
e_{j,-i}^{}:=[e_{j,-k}^{},\,e_{k,-i}^{}]_{q}~, 
\label{A9}
\end{equation}
where $1\le i<k<j\le l$. It should be stressed that the structure of the
composite root vectors is independent of choice of the index $k$ in the
r.h.s. of the definition (\ref{A9}). In particular, we have
\begin{equation}
\begin{array}{rcccl}
e_{i,-j}^{}&\!\!:=\!\!\!&[e_{i,-i-1}^{},e_{i+1,-j}^{}]_{q^{\!-1}}&\!\!\!=
\!\!\!&[e_{i,-j+1}^{},\,e_{j-1,-j}^{}]_{q^{-1}},
\\[7pt]
e_{j,-i}^{}&\!\!:=\!\!&[\,e_{j,-i-1}^{},\,e_{i+1,-i}^{}\,]_{q}&\!\!=
\!\!&[e_{j,-j+1}^{},\,e_{j-1,-i}^{}]_{q},
\end{array}
\label{A10}
\end{equation}
where $2\le i\!+\!1<j\le l$. General properties of the Cartan-Weyl basis
$\{e_{i,-j}^{}\}$ can be found in \cite{KT2}.

As it was noted above the Dynkin diagrams of the non-twisted affine
algebras can be also used for classification of the Drinfeldians and
the Yangians. In the case of $\mathfrak{sl}_{l}^{}$ ($\mathfrak{gl}_{l}^{}$),
the Dynkin diagram of the corresponding affine Lie algebra
$A^{(1)}_{l-1}\simeq\widehat{\mathfrak{sl}}_{l}^{}$
($\widehat{\mathfrak{gl}}_{l}^{}$) is presented by the picture \cite{Kac}:
\\
\begin{eqnarray}
\mbox{\begin{picture}(100,50)
\put(-40,0){\circle{5}}
\put(-37,0){\line(1,0){44}}
\put(-43,-15){\fs$\alpha_{1}^{}$}
\put(10,0){\circle{5}}
\put(13,0){\line(1,0){35}}
\put(7,-15){\fs$\alpha_{2}^{}$}
\put(53,0){$\ldots$}
\put(107,0){\line(-1,0){35}}
\put(110,0){\circle{5}}
\put(113,0){\line(1,0){44}}
\put(100,-15){\fs$\alpha_{l-2}^{}$}
\put(160,0){\circle{5}}
\put(157,-15){\fs$\alpha_{l-1}^{}$}
\put(60,49){\circle{5}}
\put(-38,2){\line(2,1){95}}
\put(158,2){\line(-2,1){95}}
\put(50,60){\fs$\delta\!-\!\theta$}
\end{picture}}
\label{A12'}
\nonumber
\end{eqnarray}
\nopagebreak
\\
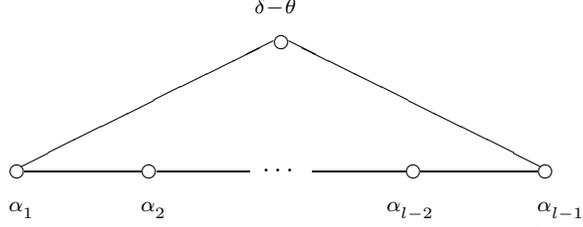
\centerline{\fs Fig. 3. The Dynkin diagram of the affine Lie algebra
$\widehat{\mathfrak{sl}}_{l}$.}
\\

The defining relations for the generators of  $D_{q\eta}(\mathfrak{g})$
presented in the previous section depend explicitly on the choice of the
element $\tilde{e}_{-\theta}\in U_{q}(\mathfrak{g})$ of the weight $-\theta$,
such that $\mathop{\lim}\limits_{q\rightarrow1}\tilde{e}_{-\theta}=
e_{-\theta}'\in\,\mathfrak{g}$.
Here we present specification to the case
$\mathfrak{g}=\mathfrak{gl}_{l}^{}$ and we set
\begin{equation}
\tilde{e}_{-\theta}=q^{e_{1,-1}\!+e_{l,-l}}e_{l,-1}~.
\label{A11}
\end{equation}
After some calculations we obtain the following result.
\begin{theorem}
The Drinfeldian $D_{q\eta}(\mathfrak{gl}_{l}^{})$ ($l\ge3$) is generated
(as a unital associative algebra over $\mathbb C[[\log q,\eta]]$) by the
algebra $U_{q}(\mathfrak{gl}_{l})$ and the elements $\xi_{\delta-\theta}^{}$,
$q^{\pm\hat{c}}:=q^{\pm h_\delta}$ with the relations:
\begin{eqnarray}
[q^{\pm\hat{c}},\,{\rm everything}]\!\!&=\!\!&0~,
\label{A12}
\\[5pt]
q^{\pm e_{1,-1}}\xi_{\delta-\theta}^{}\!\!&=\!\!&q^{\mp1}
\xi_{\delta-\theta}^{}q^{\pm e_{1,-1}}~,
\label{A13}
\\[5pt]
q^{\pm e_{i,-i}}\xi_{\delta-\theta}^{}\!\!&=\!\!&\xi_{\delta-\theta}^{}\,
q^{\pm e_{i,-i}}\quad{\rm for}\; i=2,3,\ldots,l\!-1,
\label{A14}
\\[5pt]
q^{\pm e_{l,-l}}\xi_{\delta-\theta}^{}\!\!&=\!\!&q^{\pm1}
\xi_{\delta-\theta}^{}\,q^{\pm e_{l,-l}}~,
\label{A15}
\\[5pt]
[\xi_{\delta-\theta}^{},\,e_{i+1,-i}^{}]\!\!&=\!\!&0
\qquad{\rm for}\; i=1,2,\ldots,l-1,
\label{A16}
\\[5pt]
[e_{i,-i-1}^{},\,\xi_{\delta-\theta}^{}]\!\!&=\!\!&0
\qquad{\rm for}\; i=2,3,\ldots,l-2,
\label{A17}
\\[5pt]
[e_{1,-2}^{},\,[e_{1,-2}^{},\,\xi_{\delta-\theta}^{}]_{q}]_{q}\!\!&=\!\!&0~,
\label{A18}
\\[5pt]
[[\xi_{\delta-\theta}^{},\,e_{l-1,-l}^{}]_{q},e_{l-1,-l}^{}]_{q}\!\!&=\!\!&0~,
\label{A19}
\end{eqnarray}
\begin{equation}
\begin{array}{rcl}
[[e_{1,-2}^{},\,\xi_{\delta-\theta}^{}]_{q},\xi_{\delta-\theta}^{}]_{q}\!\!&=
\!\!&\eta\,q^{e_{1,-1}+e_{l,-l}}\,\times
\\[5pt]
&&\!\!\!\!\!\!\!\!\!\!\!\!\!\!\times
\left(q^{-2}[e_{1,-2}^{},\,e_{l,-1}^{}]\xi_{\delta-\theta}^{}
 -e_{l,-1}^{}[e_{1,-2}^{},\,\xi_{\delta-\theta}^{}]_{q}\right),
\label{A20}
\end{array}
\end{equation}\
\begin{equation}
\begin{array}{rcl}
[[\xi_{\delta-\theta}^{},\,[\xi_{\delta-\theta}^{},\,
e_{l-1,-l}^{}]_{q}]_{q}\!\!&=\!\!&\eta\,q^{e_{1,-1}^{}+e_{l,-l}^{}+1}
\,\times
\\[5pt]
&&\!\!\!\!\!\!\!\!\!\!\!\!\!\times
\left(q\,[e_{l,-1}^{},\,e_{l-1,-l}^{}]\xi_{\delta-\theta}^{}
-e_{l,-1}^{}[\xi_{\delta-\theta}^{},e_{l-1,-l}^{}]_{q}\right).
\label{A21}
\end{array}
\end{equation}
The Hopf structure of $D_{q\eta}(\mathfrak{gl}_{l}^{})$ is defined by the
formulas (\ref{A5})--(\ref{A7}) for $U_{q}(\mathfrak{gl}_{l}^{})$ (i.e.
$\Delta_{q\eta}(x)=\Delta_{q}(x)$, $S_{q\eta}(x)=S_{q}(x)$ for
$(x\in U_{q}(gl_{l}^{})$), and $\Delta_{q}(q^{\pm\hat{c}})=
q^{\pm\hat{c}}\otimes q^{\pm\hat{c}}$, $S_{q}(q^{\pm\hat{c}})=q^{\mp\hat{c}}$.
The comultiplication and the antipode of $\xi_{\delta-\theta}^{}$ are given
by
\begin{equation}
\begin{array}{rcl}
\!\!\!\!\!\Delta_{q\eta}(\xi_{\delta-\theta}^{})\!\!&=\!\!&
\xi_{\delta-\theta}^{}\otimes1+
q^{e_{1,-1}-e_{l,-l}-\hat{c}}\otimes\xi_{\delta-\theta}^{}
\\[4pt]
&&+\eta\Bigl(e_{l,-1}^{}q^{e_{l,-l}}\otimes[e_{1,-1}^{}]+
\Bigl[\mbox{\ls$\frac{1}{2}$}\hat{c}+e_{l,-l}^{}\Bigr]\,
q^{-\frac{1}{2}\hat{c}}\otimes e_{l,-1}^{}q^{e_{l,-l}}
\\[4pt]
&&+{\displaystyle\sum\limits_{i=2}^{l-1}}e_{l,-i}^{}q^{e_{l,-l}}\otimes
e_{i,-1}^{}q^{e_{i,-i}}\Bigr)\Bigl(q^{e_{1,-1}}\otimes q^{e_{1,-1}}\Bigr)~,
\label{A22}
\end{array}
\end{equation}
\begin{equation}
\begin{array}{rcl}
&&\!\!\!\!\!\!\!\!\!\!\!\!\!
S_{q\eta}(\xi_{\delta-\theta}^{})\,=
-q^{\hat{c}-e_{1,-1}+e_{l,-l}}\xi_{\delta-\theta}^{}
\\[4pt]
&&+\,\eta\left(\Bigl[\mbox{\ls$\frac{1}{2}$}\hat{c}+
e_{1,-1}^{}+e_{l,-l}^{}+1\Bigl]
q^{\frac{1}{2}\hat{c}-e_{1,-1}\!+e_{l,-l}-1}e_{l,-1}^{}\right.
\\[4pt]
&&+\,{\displaystyle\sum\limits_{k=1}^{l-1}}\;q^{-k}(q-q^{-1})^{k-1}\!\!\!\!\!
{\displaystyle\sum\limits_{l-1\ge
i_{k}>\ldots>i_{1}\ge2}}\!\!\!\!\!
e_{l,-i_{k}}^{}e_{i_{k},-i_{k-1}}^{}\cdots e_{i_{1},-1}^{}q^{-2e_{1,-1}}\Bigr),
\end{array}
\label{A23}
\end{equation}
where we use the notation $[x]:=(q^{x}-q^{-x})/(q-q^{-1})$.
\label{DYP1}
\end{theorem}
It is not difficult to check that the substitution
$\xi_{\delta-\theta}^{}=q^{e_{1,-1}+e_{l,-l}}e_{l,-1}^{}$
satisfies the relations (\ref{A12})--(\ref{A21}), i.e. there is
a simple homomorphism
$D_{q\eta}(\mathfrak{gl}_{l}^{})\to U_{q}(\mathfrak{gl}_{l}^{})$.
Moreover, the both sides of the relations (\ref{A20}) and (\ref{A21})
are equal to zero independently. Therefore, we can construct a "evaluation
representation" $\rho_{ev}$ of $D_{q\eta}(\mathfrak{gl}_{l}^{})$ in
$U_q(\mathfrak{gl}_{l}^{})\otimes \mathbb C[u]$ as follows
\begin{equation}
\begin{array}{rcccl}
\rho_{ev}(q^{\hat{c}})\!\!&=\!\!&1~,\qquad\quad
\rho_{ev}(\xi_{\delta-\theta}^{})\!\!&=\!\!&u\,
q^{e_{1,-1}+e_{l,-l}}e_{l,-1}^{}~,
\\[5pt]
\rho_{ev}(q^{\pm e_{i,-i}})\!\!&=\!\!&q^{\pm e_{i,-i}},\quad
\rho_{ev}(e_{i,-j}^{})\!\!&=\!\!&e_{i,-j}^{}\quad\;\; (|i-j|=1)~.
\end{array}
\label{A24}
\end{equation}

It is obvious that
\begin{equation}
D_{q\eta=0}(\mathfrak{gl}_{l}^{})\simeq U_{q}(\mathfrak{gl}_{l}^{}[u])
\label{A25}
\end{equation}
as Hopf algebras. If $q\to $1 then the limit Hopf algebra
$D_{q=1\eta}(\mathfrak{gl}_{l}^{})$ is isomorphic to the Yangian
$Y_{\eta}(\mathfrak{gl}_{l}^{})$ with $\hat{c}\neq 0$) \cite{T2}:
\begin{equation}
D_{q=1\eta}(\mathfrak{gl}_{l}^{})\simeq Y_{\eta}(\mathfrak{gl}_{l}^{})~.
\label{A26}
\end{equation}
By setting $q=1$ in (\ref{A12})--(\ref{A23}), we obtain the defining
relations of the Yangian $Y_{\eta}(\mathfrak{gl}_{l}^{})$ and its Hopf
structure in the Chevalley basis. This result can be formulated as the
theorem.
\begin{theorem}
The Yangian $Y_{\eta}(\mathfrak{gl}_{l}^{})$ ($l\ge3$) is generated
(as a unital associative algebra over $\mathbb C[\eta]$) by the algebra
$U(\mathfrak{gl}_{l}^{})$ and the elements $\xi_{\delta-\theta}^{}$,
$\hat{c}$ with the relations:
\begin{eqnarray}
[\hat{c},\,{\rm everything}]\!\!&=\!\!&0~,
\label{A27}
\\[4pt]
[e_{1,-1}^{},\,\xi_{\delta-\theta}^{}]\!\!&=\!\!&-\xi_{\delta-\theta}^{}~,
\label{A28}
\\[4pt]
[e_{l,-l}^{},\,\xi_{\delta-\theta}^{}]\!\!&=\!\!&\xi_{\delta-\theta}^{}~,
\label{A29}
\\[4pt]
[e_{i,-i}^{},\,\xi_{\delta-\theta}^{}]\!\!&=\!\!&0
\quad\;{\rm for}\; i=2,3,\ldots,l-1,
\label{A30}
\\[4pt]
[\xi_{\delta-\theta}^{},\,e_{i+1,-i}^{}]\!\!&=\!\!&0
\quad\;{\rm for}\; i=1,2,\ldots,l-1,
\label{A31}
\\[4pt]
[e_{i,-i-1}^{},\,\xi_{\delta-\theta}^{}]\!\!&=\!\!&0
\quad\;{\rm for}\; i=2,3,\ldots,l-2,
\label{A32}
\\[4pt]
[e_{1,-2}^{},\,[e_{1,-2},\,\xi_{\delta-\theta}^{}]]\!\!&=\!\!&0~,
\label{A33}
\\[4pt]
[[\xi_{\delta-\theta}^{},\,e_{l-1,-l}^{}],\,e_{l-1,-l}^{}]\!\!&=\!\!&0~,
\label{A34}
\end{eqnarray}
\begin{eqnarray}
[[e_{1,-2}^{},\,\xi_{\delta-\theta}^{}],\,\xi_{\delta-\theta}^{}]\!\!&=\!\!&
\eta\Bigl([e_{1,-2}^{},\,e_{l,-1}^{}]\,\xi_{\delta-\theta}^{}-
e_{l,-1}^{}\,[e_{1,-2}^{},\,\xi_{\delta-\theta}^{}]\Bigr),
\label{A35}
\\[4pt]
\phantom{aaaaa}
[[\xi_{\delta-\theta}^{},\,[\xi_{\delta-\theta}^{},\,e_{l-1,-l}^{}]]\!\!&=\!\!&
\eta\Bigl([e_{l,-1}^{},\,e_{l-1,-l}^{}]\,\xi_{\delta-\theta}^{}-
e_{l,-1}^{}\,[\xi_{\delta-\theta}^{},\,e_{l-1,-l}^{}]\Bigr).
\label{A36}
\end{eqnarray}
The Hopf structure of the Yangian is trivial for $U(\mathfrak{gl}_{l}^{})\oplus
\mathbb{C}\,\hat{c}\subset Y_{\eta}(\mathfrak{gl}_{l}^{})$
(i.e. $\Delta_{\eta}(x)=x\otimes 1+1\otimes x$, $S_{\eta}(x)=-x$ for
$x\in\mathfrak{gl}_{l}^{}\oplus\mathbb{C}\,\hat{c}$) and it is not trivial
for the element $\xi_{\delta-\theta}^{}$:
\begin{eqnarray}
\phantom{aaaaaaa}
\Delta_{\eta}(\xi_{\delta-\theta}^{})\!\!&=\!\!&
\xi_{\delta-\theta}^{}\otimes 1+1\otimes\xi_{\delta-\theta}^{}
+\eta\,\Big(\mbox{\Large$\frac{1}{2}$}\,\hat{c}\otimes e_{l,-1}^{}
+\sum\limits_{i=1}^{l}e_{l,-i}^{}\otimes e_{i,-1}^{}\Big),
\label{A37}
\\[4pt]
S_{\eta}(\xi_{\delta-\theta}^{})\!\!&=\!\!&-\xi_{\delta-\theta}^{}
+\eta\Big(\mbox{\Large$\frac{1}{2}$}\,\hat{c}\,e_{l,-1}^{}
+\sum\limits_{i=1}^{l}e_{l,-i}^{}\,e_{i,-1}^{}\Big).
\label{A38}
\end{eqnarray}
\label{AP2}
\end{theorem}

\setcounter{equation}{0}
\section{Drinfeldians and Yangians for Lie algebras $B_l^{}$ ($l\geq3$)}
The Dynkin diagram of the simple Lie algebra
$B_l\simeq\mathfrak{so}_{2l+1}^{}$ is presented on the picture:
\\
\begin{eqnarray}
\mbox{\begin{picture}(100,10)
\put(-40,0){\circle{5}}
\put(-37,0){\line(1,0){44}}
\put(-43,-15){\fs$\alpha_{1}^{}$}
\put(10,0){\circle{5}}
\put(13,0){\line(1,0){35}}
\put(7,-15){\fs$\alpha_{2}^{}$}
\put(55,0){$\ldots$}
\put(107,0){\line(-1,0){35}}
\put(110,0){\circle{5}}
\put(111,2){\line(1,0){46}}
\put(111,-2){\line(1,0){46}}
\put(154,-3){$\bigr>$}
\put(100,-15){\fs$\alpha_{l-1}^{}$}
\put(160,0){\circle{5}}
\put(157,-15){\fs$\alpha_{l}^{}$}
\end{picture}}
\label{B1'}
\nonumber
\end{eqnarray}
\nopagebreak
\\[2pt]
\centerline{\fs Fig. 4. The Dynkin diagram of the Lie algebra
$\mathfrak{so}_{2l+1}^{}$.}
\\

In the terms of the orthonormalized basis $\epsilon_i^{}$ ($i=1,2,\ldots,l$)
the root systems $\Pi$, $\Delta_{+}^{}$ and $\Delta$ of
$\mathfrak{so}_{2l+1}^{}$ are given as follows:
\begin{eqnarray}
\Pi\!\!&=\!\!&\{\alpha_1^{}\!=\epsilon_1^{}\!-\epsilon_2^{},\,
\alpha_2^{}\!=\epsilon_2^{}\!-\epsilon_3^{},\,\ldots,\alpha_{l-1}^{}\!=
\epsilon_{l-1}^{}\!-\epsilon_{l}^{},\,\alpha_l^{}\!=\epsilon_{l}^{}\},
\label{B1}
\\[4pt]
\Delta_{+}\!\!\!\!&=\!\!&\{\epsilon_i^{}\pm\epsilon_j^{},\,
\epsilon_k^{}\,|\;1\le i<j\le l;\, k=1,2,\ldots,l\},
\label{B2}
\\[4pt]
\Delta_{}\!\!&=\!\!&\Delta_{+}\cup\,(-\Delta_{+})\;=\;
\{\pm\epsilon_i\pm\epsilon_j,\,\pm\epsilon_k^{}\,|\,i\neq j;
\,i,j,k=1,2,\ldots,l\}.
\label{B3}
\end{eqnarray}
The root $\theta:=\alpha_1^{}+2\alpha_2^{}+\ldots+2\alpha_{l}^{}=
\epsilon_{1}^{}+\epsilon_{2}^{}$ is maximal.

The quantum algebra $U_{q}(\mathfrak{so}_{2l+1}^{}\!)$ is generated by the
Chevalley $U_q(gl_l^{})$-elements $e_{\!i,-i-1}^{}\!\!~:=~\!\!e_{\epsilon_i-
\epsilon_{i+1}}^{}$, $e_{i+1,-i}^{}:=e_{\epsilon_{i+1}-\epsilon_{i}}^{}$
$(i=1,2,\ldots,l-1)$, $q^{\pm e_{i,-i}}$ $(i=1,2,\ldots,l)$
with the relations (\ref{A4}), and the additional elements
$e_{\pm l}^{}:=e_{\pm\epsilon_{l}^{}}^{}$ with the relations:
\begin{equation}
\begin{array}{rcl}
q^{e_{i,-i}}e_{\pm l}^{}\,q^{-e_{i,-i}}\!\!&=\!\!&
q^{\pm\delta_{i,l}}e_{\pm l}^{}~,
\\[5pt]
[e_{l}^{},\,e_{-l}^{}]\!\!&=\!\!&
\mbox{\ls$\frac{q^{e_{l,-l}^{}}\,-\,q^{-e_{l,-l}^{}}}{q\,-\,q^{-1}}$}~,
\\[5pt]
[e_{i,-j}^{},e_{l}^{}]\!\!&=\!\!&0\qquad{\rm for}\;\,j\neq l,\;|i-j|=1~,
\\[5pt]
[e_{i,-j}^{},e_{-l}^{}]\!\!&=\!\!&0\qquad{\rm for}\;\,i\neq l,\;|i-j|=1~,
\\[5pt]
[[[e_{l-1,-l}^{},\,e_{l}^{}]_{q}^{},\,e_{l}^{}]_{q}^{},\,
e_{l}^{}]_{q}^{}\,\!\!&=\!\!&0~,
\\[5pt]
[[[e_{l,-l+1}^{},\,e_{-l}^{}]_{q}^{},\,e_{-l}]_{q}^{},\,
e_{-l}^{}]_{q}^{}\!\!&=\!\!&0~.
\label{B4}
\end{array}
\end{equation}
The Hopf structure on $U_{q}(\mathfrak{so}_{2l+1}^{})$ is given by the formulas
(\ref{A5})--(\ref{A7}) and also
\begin{equation}
\begin{array}{rcl}
\Delta_{q}(e_{l}^{})\!\!&=\!\!&e_{l}^{}\otimes1+q^{-e_{l,-l}^{}}\otimes e_{l}^{}~,
\\[5pt]
\Delta_{q}(e_{-l}^{})\!\!&=\!\!&e_{-l}^{}\otimes q^{e_{l,-l}^{}}+
1\otimes e_{-l}^{}~,
\end{array}
\label{B5}
\end{equation}
\begin{equation}
\begin{array}{rcccl}
S_{q}(e_{l}^{})\!\!&=\!\!&-q^{e_{l,-l}^{}}\,e_{l}^{},\qquad
S_{q}(e_{-l}^{})\!\!&=\!\!&-e_{-l}^{}\,q^{-e_{l,-l}^{}},
\end{array}
\label{B6}
\end{equation}
\begin{equation}
\varepsilon_{q}(e_{\pm l}^{})=0~.\phantom{aaaaaaaa}
\end{equation}
For construction of the composite root vectors $e_{i-j}^{}$ ($|i-j|>1$),
$e_{i,j}$, $e_{-j,-i}$ ($i<j$), and $e_{\pm i}$ $(1\leq i\leq l-1)$
we fix the following normal ordering of the positive root system:
$\Delta_{+}^{}$ 
\begin{equation}
\begin{array}{lcr}
&(\epsilon_1^{}-\epsilon_2^{},\epsilon_1^{}-\epsilon_3^{},\ldots,
\epsilon_1^{}-\epsilon_l^{},\epsilon_1^{},\epsilon_1^{}+\epsilon_l^{},
\ldots,\epsilon_{1}^{}+\epsilon_{3}^{},
\epsilon_{1}^{}+\epsilon_{2}^{}),(\epsilon_2^{}-\epsilon_3^{},&
\\[5pt]
&\epsilon_2^{}-\epsilon_4^{},\ldots,
\epsilon_2^{}-\epsilon_l^{},\epsilon_2^{},\epsilon_2^{}+\epsilon_l^{},
\ldots,\epsilon_2^{}+\epsilon_4^{},\epsilon_{2}^{}+\epsilon_{3}^{}),
\ldots,(\epsilon_{l-2}^{}-\epsilon_{l-1}^{},&
\\[5pt]
&\epsilon_{l-2}^{}-
\epsilon_{l}^{},\epsilon_{l-2}^{},\epsilon_{l-2}^{}+\epsilon_{l}^{},
\epsilon_{l-2}^{}+\epsilon_{l-1}^{}),(\epsilon_{l-1}^{}-\epsilon_{l}^{},
\epsilon_{l-1}^{},\epsilon_{l-1}^{}+\epsilon_{l}^{}),\,\epsilon_l^{}.&
\end{array}
\label{B8}
\end{equation}
According to this ordering we set
\begin{equation}
e_{i,-j}^{}\;:=\;[e_{i,-i-1}^{},\,e_{i+1,-j}^{}]_{q^{-1}},
\qquad e_{j,-i}^{}\;:=\;[e_{j,-i-1}^{},\,e_{i+1,-i}^{}]_{q}
\label{B9a}
\end{equation}
for $2\le i+1<j\le l$,
\begin{equation}
e_{i}^{}\;:=\;[e_{i,-l}^{},\,e_{l}^{}]_{q^{-1}},
\qquad e_{-i}^{}\;:=\;[e_{-l}^{},\,e_{l,-i}^{}]_{q}
\label{B9b}
\end{equation}
for $1\le i<l$,
\begin{equation}
e_{i,l}^{}\;:=\;[e_{i}^{},\,e_{l}^{}]_{q^{-1}},
\qquad e_{-l,-i}^{}\;:=\;[e_{-l}^{},\,e_{-i}^{}]_{q}
\label{B9c}
\end{equation}
for $1\le i<l$,
\begin{equation}
e_{i,j}^{}\;:=\;[e_{i,j+1}^{},\,e_{j,-j-1}^{}]_{q^{-1}},
\qquad e_{-j,-i}^{}\;:=\;[e_{j+1,-j}^{},\,e_{-j-1,-i}^{}]_{q}
\label{B9}
\end{equation}
for $1\le i<j<l$.

The extended Dynkin diagram of $B_l^{}$ (or the Dynkin diagram
of the corresponding affine Lie algebra $\widehat{\mathfrak{so}}_{2l+1}^{}$)
is presented by the picture: 
\\
\begin{eqnarray}
\mbox{\begin{picture}(125,45)
\put(-50,0){\circle{5}}
\put(-47,0){\line(1,0){44}}
\put(-53,-15){\fs$\alpha_{1}^{}$}
\put(0,0){\circle{5}}
\put(3,0){\line(1,0){44}}
\put(-3,-15){\fs$\alpha_{2}^{}$}
\put(0,40){\circle{5}}
\put(0,3){\line(0,1){34}}
\put(10,38){\fs$\delta-\theta$}
\put(50,0){\circle{5}}
\put(53,0){\line(1,0){35}}
\put(47,-15){\fs$\alpha_{3}^{}$}
\put(93,0){$\ldots$}
\put(147,0){\line(-1,0){35}}
\put(150,0){\circle{5}}
\put(151,2){\line(1,0){46}}
\put(151,-2){\line(1,0){46}}
\put(194,-3){$\bigr>$}
\put(140,-15){\fs$\alpha_{l-1}^{}$}
\put(200,0){\circle{5}}
\put(197,-15){\fs$\alpha_{l}^{}$}
\end{picture}}
\label{B11'}
\nonumber
\end{eqnarray}
\nopagebreak
\\
\centerline{\fs Fig. 5. The Dynkin diagram of the Lie algebra
$\widehat{\mathfrak{so}}_{2l+1}^{}$.}
\\

Specializing the general formulas (\ref{tolsD5})--(\ref{tolsD8}) to
the case $\mathfrak{g}=\mathfrak{so}_{2l+1}^{}$ with the vector
\begin{equation}
\tilde{e}_{-\theta}^{}=q^{e_{1,-1}^{}\!+e_{2,-2}^{}}e_{-2,-1}~,
\label{B10}
\end{equation}
after some calculations we obtain the following result.
\begin{theorem}
The Drinfeldian $D_{q\eta}(\mathfrak{so}_{2l+1}^{})$ ($l\ge3$) is generated (as
a unital associative algebra over $\mathbb C[[\log q,\eta]]$) by the algebra
$U_{q}(\mathfrak{so}_{2l+1}^{})$ and the elements $\xi_{\delta-\theta}^{}$,
$q^{\pm\hat{c}}:=q^{\pm h_\delta}$ with the relations:
\begin{eqnarray}
[q^{\pm\hat{c}},\,{\rm everything}]\!\!&=\!\!&0~,
\label{B11}
\\[5pt]
q^{\pm e_{i,-i}}\xi_{\delta-\theta}^{}\!\!&=\!\!&q^{\mp1}
\xi_{\delta-\theta}^{}q^{\pm e_{i,-i}}\quad {\rm for}\;\,i=1,2,
\label{B12}
\\[5pt]
q^{\pm e_{i,-i}}\xi_{\delta-\theta}^{}\!\!&=\!\!&\xi_{\delta-\theta}^{}\,
q^{\pm e_{i,-i}}\qquad\;{\rm for}\;\, i=3,4,\ldots,l,
\label{B13}
\end{eqnarray}
\begin{equation}
\begin{array}{rcl}
[e_{2,-1}^{},\,\xi_{\delta-\theta}^{}]\!\!&=\!\!&
\eta q^{e_{1,-1}^{}+e_{2,-2}^{}}\,\times
\\[3pt]
&&\!\!\!\!\!\!\times\,
\Bigl(\sum\limits_{k=3}^{l}(-1)^kq^{k-3}e_{-k,-1}^{}e_{k,-1}^{}
\!+\mbox{\ls$\frac{(-1)^{l-1}q^{l-1}}{q+1}$}e_{-1}^2\Bigr),
\label{B14}
\end{array}
\end{equation}
\begin{equation}
\begin{array}{rcl}
[e_{1,-2}^{},\,\xi_{\delta-\theta}^{}]\!\!&=\!\!&
\eta q^{2e_{2,-2}^{}}\,\times
\\[3pt]
&&\!\!\!\!\!\!\!
\times\,\Bigl(\sum\limits_{k=3}^{l}(-1)^kq^{k-3}e_{-k,-2}^{}e_{k,-2}^{}
\!+\mbox{\ls$\frac{(-1)^{l-1}q^{l-1}}{q+1}$}e_{-2}^2\Bigr),
\label{B15}
\end{array}
\end{equation}
\begin{eqnarray}
[e_{3,-2}^{},\,\xi_{\delta-\theta}^{}]\!\!&=\!\!&0~,
\\[5pt]
\label{B16}
[e_{2,-3}^{},\,[e_{2,-3}^{},\,\xi_{\delta-\theta}^{}]_{q}]_{q}\!\!&=\!\!&0~,
\label{B17}
\\[9pt]
[[e_{2,-3}^{},\,\xi_{\delta-\theta}^{}]_{q},\,
\xi_{\delta-\theta}^{}]_{q}\!\!&=\!\!&0~,
\label{B18}
\\[5pt]
[e_{i,-j}^{},\,\xi_{\delta-\theta}^{}]\!\!&=\!\!&0
\quad (3\leq i,j\leq l,\;|i-j|=1),
\label{B19}
\\[5pt]
[e_{\pm l}^{},\,\xi_{\delta-\theta}^{}]\!\!&=\!\!&0~.
\label{B20}
\end{eqnarray}
\end{theorem}
Since the explicit formulas of $\Delta_{q\eta}(\xi_{\delta-\theta}^{})$ and
$S_{q\eta}(\xi_{\delta-\theta}^{})$ for the cases $B_l$, $C_l$ and $D_l$
are cumbersome they are not written here. Analogous formulas are also not
given for the corresponding Yangians.

By setting $q=1$ in (\ref{B11})--(\ref{B20}), we obtain the defining
relations of the Yangian $Y_{\eta}(\mathfrak{so}_{2l+1}^{})$
in the Chevalley basis. This result is formulated as the theorem.
\begin{theorem}
The Yangian $Y_{\eta}(\mathfrak{so}_{2l+1}^{})$ ($l\ge 3$) is generated
(as a unital associative algebra over $\mathbb C[\eta]$) by the algebra
$U(\mathfrak{so}_{2l+1}^{})$ and the elements $\xi_{\delta-\theta}^{}$,
$\hat{c}$ with the relations:
\begin{eqnarray}
[\hat{c},\,{\rm everything}]\!\!&=\!\!&0~,
\label{B21}
\\[5pt]
[e_{i,-i},\,\xi_{\delta-\theta}^{}]\!\!&=\!\!&\xi_{\delta-\theta}^{}
\quad {\rm for}\;\, i=1,2,
\label{B22}
\\[5pt]
[e_{i,-i}^{}\xi_{\delta-\theta}^{}]\!\!&=\!\!&0
\qquad\;{\rm for}\;\,i=3,4,\ldots,l,
\label{B23}
\\[1pt]
[e_{2,-1}^{},\,\xi_{\delta-\theta}^{}]\!\!&=\!\!&
\eta\Bigl(\sum_{k=3}^{l}(-1)^ke_{-k,-1}^{}e_{k,-1}^{}
+\mbox{\ls$\frac{(-1)^{l-1}}{2}$}\,e_{-1}^2\Bigr),
\label{B24}
\\[1pt]
[e_{1,-2}^{},\,\xi_{\delta-\theta}^{}]\!\!&=\!\!&
\eta\Bigl(\sum_{k=3}^{l}(-1)^ke_{-k,-2}^{}e_{k,-2}^{}
+\mbox{\ls$\frac{(-1)^{l-1}}{2}$}\,e_{-2}^2\Bigr),
\label{B25}
\\[1pt]
[e_{3,-2}^{},\,\xi_{\delta-\theta}^{}]\!\!&=\!\!&0~,
\\[5pt]
\label{B26}
[e_{2,-3}^{},\,[e_{2,-3}^{},\,\xi_{\delta-\theta}^{}]]\!\!&=\!\!&0~,
\label{B27}
\\[5pt]
[[e_{2,-3}^{},\,\xi_{\delta-\theta}^{}],\,\xi_{\delta-\theta}^{}]\!\!&=\!\!&0~,
\label{B28}
\\[5pt]
[e_{i,-j}^{},\,\xi_{\delta-\theta}^{}]\!\!&=\!\!&0
\qquad
(3\leq i,j\leq l,\;|i-j|=1),
\label{B29}
\\[5pt]
[e_{\pm l}^{},\,\xi_{\delta-\theta}^{}]\!\!&=\!\!&0~.
\label{B30}
\end{eqnarray}
\end{theorem}

\setcounter{equation}{0}
\section{Drinfeldians and Yangians for Lie algebras $C_l^{}$ ($l\geq2$)}
The Dynkin diagram of the simple Lie algebra
$C_l\simeq\mathfrak{sp}_{2l}^{}:=\mathfrak{sp}(2l,\mathbb C)$
is presented on the picture:
\\
\begin{eqnarray}
\mbox{\begin{picture}(105,10)
\put(-40,0){\circle{5}}
\put(-37,0){\line(1,0){44}}
\put(-43,-15){\fs$\alpha_{1}^{}$}
\put(10,0){\circle{5}}
\put(13,0){\line(1,0){35}}
\put(7,-15){\fs$\alpha_{2}^{}$}
\put(55,0){$\ldots$}
\put(107,0){\line(-1,0){35}}
\put(110,0){\circle{5}}
\put(113,2){\line(1,0){46}}
\put(113,-2){\line(1,0){46}}
\put(111,-3){$\bigl<$}
\put(100,-15){\fs$\alpha_{l-1}^{}$}
\put(160,0){\circle{5}}
\put(157,-15){\fs$\alpha_{l}^{}$}
\end{picture}}
\label{C1'}
\nonumber
\end{eqnarray}
\nopagebreak
\\
\nopagebreak
\centerline{\fs Fig. 6. The Dynkin diagram of the Lie algebra
$\mathfrak{sp}_{2l}^{}$.}
\\

In the terms of the orthonormalized basis $\epsilon_i^{}$
($i=1,2,\ldots,l$) the root systems $\Pi$, $\Delta_{+}^{}$ and $\Delta$ of
$\mathfrak{sp}_{2l}^{}$ are presented as follows:
\begin{eqnarray}
\Pi\!\!&=\!\!&\{\alpha_1^{}\!=\epsilon_1^{}\!-\epsilon_2^{},\,
\alpha_2^{}\!=\epsilon_2^{}\!-\epsilon_3^{},\,\ldots,\alpha_{l-1}^{}\!=
\epsilon_{l-1}^{}\!-\epsilon_{l}^{},\,\alpha_l^{}\!=2\epsilon_{l}^{}\},
\label{C1}
\\[4pt]
\Delta_{+}\!\!\!\!&=\!\!&\{\epsilon_i^{}\pm\epsilon_j^{},\,
2\epsilon_k^{}\,|\;1\le i<j\le l;\, k=1,2,\ldots,l\},
\label{C2}
\\[4pt]
\Delta_{}\!\!&=\!\!&\Delta_{+}\cup\,(-\Delta_{+})\;=\;
\{\pm\epsilon_i\pm\epsilon_j,\,\pm 2\epsilon_k^{}\,|\,
i\neq j;\,i,j,k=1,2,\ldots,l\}.
\label{C3}
\end{eqnarray}
The root $\theta:=2\alpha_1^{}+\ldots+2\alpha_{l-1}^{}+
\alpha_{l}^{}=2\epsilon_{1}^{}$ is maximal.

The quantum algebra $U_{q}(\mathfrak{sp}_{2l}^{})$ is generated by the
Chevalley $U_q(gl_l^{})$-elements $e_{i,-i-1}^{}:=e_{\epsilon_i-
\epsilon_{i+1}}^{}$, $e_{i+1,-i}^{}:=e_{\epsilon_{i+1}-\epsilon_{i}}^{}$
$(i=1,2,\ldots,l-1)$, $q^{\pm e_{i,-i}}$ $(i=1,2,\ldots,l)$
satisfying the relations (\ref{A4}) and the additional elements
$e_{l,l}^{}:=e_{2\epsilon_{l}^{}}^{}$, $e_{-l,-l}^{}:=
e_{-2\epsilon_{l}^{}}^{}$ with the relations:
\begin{equation}
\begin{array}{rcl}
q^{e_{i,-i}}e_{\pm l,\pm l}^{}\,q^{-e_{i,-i}}\!\!&=\!\!&
q^{\pm2\delta_{i,l}}e_{\pm l,\pm l}^{},
\\[5pt]
[e_{l,l}^{},\,e_{-l,-l}^{}]\!\!&=\!\!&
\mbox{\ls$\frac{q^{2e_{l,-l}}\,-\,q^{-2e_{l,-l}}}{q\,-\,q^{-1}}$},
\\[7pt]
[e_{i,-j}^{},e_{l,l}^{}]\!\!&=\!\!&0\quad\;{\rm for}\;\,j\neq l,\;\,|i-j|=1,
\\[5pt]
[e_{i,-j}^{},e_{-l,-l}^{}]\!\!&=\!\!&0\quad\;{\rm for}\;\,i\neq l,\;\,|i-j|=1,
\\[5pt]
[[e_{l-1,-l}^{},\,e_{l,l}^{}]_{q}^{},\,e_{l,l}^{}]_{q}\!\!&=\!\!&0~,
\\[5pt]
[[e_{l,-l+1}^{},\,e_{-l,-l}^{}]_{q}^{},\,e_{-l,-l}]_{q}\!\!&=\!\!&0~,
\end{array}
\end{equation}
\begin{equation}
\begin{array}{rcl}
[e_{l-1,-l}^{},\,[e_{l-1,-l}^{},\,[e_{l-1,-l}^{},\,
e_{l,l}^{}]_{q}]_q]_q\!\!&=\!\!&0,
\\[5pt]
[e_{l,-l+1}^{},\,[e_{l,-l+1}^{},\,[e_{l,-l+1}^{},\,
e_{-l,-l}^{}]_{q}]_q]_q\!\!&=\!\!&0.
\label{C4}
\end{array}
\end{equation}
The Hopf structure on $U_{q}(\mathfrak{sp}_{2l}^{})$ is given by the
formulas (\ref{A5})--(\ref{A7}) and also
\begin{equation}
\begin{array}{rcl}
\Delta_{q}(e_{l,l}^{})\!\!&=\!\!&e_{l,l}^{}\otimes1+
q^{-2e_{l,-l}}\otimes e_{l,l}^{},
\\[7pt]
\Delta_{q}(e_{-l,-l}^{})\!\!&=\!\!&e_{-l,-l}^{}\otimes
q^{2e_{l,-l}}+1\otimes e_{-l,-l}^{},
\end{array}
\label{C5}
\end{equation}
\begin{equation}
\begin{array}{rcccl}
S_{q}(e_{l,l}^{})\!\!&=\!\!&-q^{2e_{l,-l}}\,e_{l,l}^{},\qquad
S_{q}(e_{-l,-l}^{})\!\!&=\!\!&-e_{-l,-l}^{}\,q^{-2e_{l,-l}},
\end{array}
\label{C6}
\end{equation}
\begin{equation}
\varepsilon_{q}(e_{\pm l,\pm l}^{})=0.
\label{C7}
\end{equation}

For construction of the composite root vectors $e_{i,-j}^{}$
($|i-j|>1$), $e_{i,j}^{}$, $e_{-j,-i}^{}$ ($i\le j$)
we fix the following normal ordering of the positive root system
$\Delta_{+}^{}$: 
\begin{equation}
\begin{array}{rcl}
&&(\epsilon_1^{}\!-\epsilon_2^{},\epsilon_1^{}\!-\epsilon_3^{},
\ldots,\epsilon_1^{}\!-\epsilon_l^{},2\epsilon_1^{},
\epsilon_1^{}\!+\epsilon_l^{},\ldots,\epsilon_{1}^{}\!+\epsilon_{3}^{},
\epsilon_{1}^{}\!+\epsilon_{2}^{}),(\epsilon_2^{}\!-\epsilon_3^{},
\\[5pt]
&&\epsilon_2^{}\!-\epsilon_4^{}\ldots,\epsilon_2^{}\!-\epsilon_l^{},
2\epsilon_2^{},\epsilon_2^{}\!+\epsilon_l^{},\ldots,\epsilon_{2}^{}\!+
\epsilon_{3}^{}),\ldots,(\epsilon_{l-2}^{}\!-\epsilon_{l-1}^{},
\epsilon_{l-2}^{}\!-\epsilon_{l}^{},
\\[5pt]
&&2\epsilon_{l-2}^{},\epsilon_{l-2}^{}\!+\epsilon_{l}^{},
\epsilon_{l-2}^{}\!+\epsilon_{l-1}^{}),(\epsilon_{l-1}^{}\!-\epsilon_{l}^{},
2\epsilon_{l-1}^{},\epsilon_{l-1}^{}\!+\epsilon_{l}^{}),\,2\epsilon_l^{}.
\end{array}
\label{C8}
\end{equation}
According to this ordering we set
\begin{equation}
e_{i,-j}^{}\;:=\;[e_{i,-i-1}^{},\,e_{i+1,-j}^{}]_{q^{-1}},
\qquad e_{j,-i}^{}\;:=\;[e_{j,-i-1}^{},\,e_{i+1,-i}^{}]_{q}
\end{equation}
for $2\le i\!+\!1<j\le l$,
\begin{equation}
e_{i,l}^{}\;:=\;[e_{i,-l}^{},\,e_{l,l}^{}]_{q^{-1}},
\qquad e_{-i,-l}^{}\;:=\;[e_{-l,-l}^{},\,e_{l,-i}^{}]_{q}
\end{equation}
for $1\le i< l$, 
\begin{equation}
e_{i,i}^{}\;:=\;[e_{i,-l}^{},\,e_{i,l}^{}]_{q^{-1}},
\qquad e_{-i,-i}^{}\;:=\;[e_{-l,-i}^{},\,e_{l,-i}^{}]_{q}
\end{equation}
for $1\le i< l$,
\begin{equation}
e_{i,j}^{}\;:=\;[e_{i,j+1}^{},\,e_{j,-j-1}^{}]_{q^{-1}},
\qquad e_{-j,-i}^{}\;:=\;[e_{j+1,-j}^{},\,e_{-j-1,-i}^{}]_{q}
\label{C9}
\end{equation}
for $1\le i<j<l$.

The extended Dynkin diagram of $C_l^{}$ (or the Dynkin diagram of the
corresponding affine Lie algebra $\widehat{\mathfrak{sp}}_{2l}^{}$)
is presented by the picture: 
\\
\begin{eqnarray}
\mbox{\begin{picture}(120,10)
\put(-50,0){\circle{5}}
\put(-49,2){\line(1,0){46}}
\put(-49,-2){\line(1,0){46}}
\put(-6,-3){$\bigr>$}
\put(-60,-15){\fs$\delta-\theta$}
\put(0,0){\circle{5}}
\put(3,0){\line(1,0){44}}
\put(-3,-15){\fs$\alpha_{1}^{}$}
\put(50,0){\circle{5}}
\put(53,0){\line(1,0){35}}
\put(47,-15){\fs$\alpha_{2}^{}$}
\put(93,0){$\ldots$}
\put(147,0){\line(-1,0){35}}
\put(150,0){\circle{5}}
\put(153,2){\line(1,0){46}}
\put(153,-2){\line(1,0){46}}
\put(151,-3){$\bigl<$}
\put(140,-15){\fs$\alpha_{l-1}^{}$}
\put(200,0){\circle{5}}
\put(197,-15){\fs$\alpha_{l}^{}$}
\end{picture}}
\label{C11'}
\nonumber
\end{eqnarray}
\nopagebreak
\\[0pt]
\centerline{\fs Fig. 7. The Dynkin diagram of the Lie algebra
$\widehat{sp}_{2l}^{}$.}
\\

Specializing the general formulas and (\ref{tolsD5})--(\ref{tolsD8}) to
the case $\mathfrak{g}=\mathfrak{sp}_{2l}^{}$ with the element
\begin{equation}
\tilde{e}_{-\theta}^{}=q^{2e_{1,-1}^{}}e_{-1,-1}~,
\label{C10}
\end{equation}
after some calculations we obtain the following result.
\begin{theorem}
The Drinfeldian $D_{q\eta}(\mathfrak{sp}_{2l}^{})$ ($l\ge2$) is generated (as
a unital associative algebra over $\mathbb C[[\log q,\eta]]$) by the algebra
$U_{q}(\mathfrak{sp}_{2l}^{})$ and the elements $\xi_{\delta-\theta}^{}$,
$q^{\pm\hat{c}}:=q^{\pm h_\delta}$ with the relations:
\begin{eqnarray}
[q^{\pm\hat{c}},\,{\rm everything}]\!\!&=\!\!&0~,
\label{C11}
\\[4pt]
q^{\pm e_{1,-1}^{}}\xi_{\delta-\theta}^{}\!\!&=\!\!&q^{\mp2}
\xi_{\delta-\theta}^{}q^{\pm e_{1,-1}},
\label{C12}
\\[4pt]
q^{\pm e_{i,-i}}\xi_{\delta-\theta}^{}\!\!&=\!\!&\xi_{\delta-\theta}^{}\,
q^{\pm e_{i,-i}}\quad(i=2,3,\ldots,l),
\label{C13}
\\[4pt]
[e_{2,-1}^{},\,\xi_{\delta-\theta}^{}]\!\!&=\!\!&0~,
\label{C14}
\\[4pt]
\phantom{aaaaa}[e_{1,-2}^{},\,[e_{1,-2}^{},\,[e_{1,-2}^{},\,
\xi_{\delta-\theta}^{}]_{q}]_{q}]_{q}\!\!&=\!\!&0~,
\label{C15}
\\[4pt]
[[e_{1,-2}^{},\,\xi_{\delta-\theta}^{}]_{q},\,
\xi_{\delta-\theta}^{}]_{q}\!\!&=\!\!&0~,
\label{C16}
\\[4pt]
[e_{i,-j}^{},\,\xi_{\delta-\theta}^{}]\!\!&=\!\!&0
\quad
(2\leq i,j\leq l,\;|i-j|=1),
\label{C17}
\\[4pt]
[e_{-l,-l}^{},\,\xi_{\delta-\theta}^{}]\!\!&=\!\!&
\eta\,q^{2e_{1,-1}^{}}e_{-l,-1}^2~,
\label{C18}
\\[4pt]
[e_{l,l}^{},\,\xi_{\delta-\theta}^{}]\!\!&=\!\!&
\eta\,[2]\,q^{2(e_{1,-1}^{}+e_{l,-l}^{})-3}e_{l,-1}^2.
\label{C19}
\end{eqnarray}
\end{theorem}

By setting $q=1$ in (\ref{C11})--(\ref{C19}), we obtain the defining
relations of the Yangian $Y_{\eta}(\mathfrak{sp}_{2l}^{})$
in the Chevalley basis. This result can be formulated as the theorem.
\begin{theorem}
The Yangian $Y_{\eta}(\mathfrak{sp}_{2l}^{})$ ($l\ge2$) is generated
(as a unital associative algebra over $\mathbb C[\eta]$) by the algebra
$U(\mathfrak{sp}_{2l}^{})$ and the elements $\xi_{\delta-\theta}^{}$,
$\hat{c}$ with the relations:
\begin{eqnarray}
[\hat{c},\,{\rm everything}]\!\!&=\!\!&0~,
\label{C20}
\\[5pt]
[e_{1,-1}^{},\,\xi_{\delta-\theta}^{}]\!\!&=\!\!&-2\xi_{\delta-\theta}^{},
\label{C21}
\\[5pt]
[e_{i,-i}^{},\,\xi_{\delta-\theta}^{}]\!\!&=\!\!&0
\quad(i=2,3,\ldots,l)~,
\label{C22}
\\[5pt]
[e_{2,-1}^{},\,\xi_{\delta-\theta}^{}]\!\!&=\!\!&0~,
\label{C23}
\\[5pt]
\phantom{aaaaa}[e_{1,-2}^{},\,[e_{1,-2}^{},\,[e_{1,-2}^{},\,
\xi_{\delta-\theta}^{}]]]\!\!&=\!\!&0~,
\label{C24}
\\[5pt]
[[e_{1,-2}^{},\,\xi_{\delta-\theta}^{}],\,
\xi_{\delta-\theta}^{}]\!\!&=\!\!&0~,
\label{C25}
\\[5pt]
[e_{i,-j}^{},\,\xi_{\delta-\theta}^{}]\!\!&=\!\!&0
\quad 
(2\leq i,j\leq l,\;|i-j|=1),
\label{C26}
\\[5pt]
[e_{-l,-l}^{},\,\xi_{\delta-\theta}^{}]\!\!&=\!\!&\eta\,e_{-l,-1}^2~,
\label{C27}
\\[5pt]
[e_{l,l}^{},\,\xi_{\delta-\theta}^{}]\!\!&=\!\!&2\,\eta\,e_{l,-1}^2~.
\label{C28}
\end{eqnarray}
\end{theorem}

\setcounter{equation}{0}
\section{Drinfeldians and Yangians for Lie algebras $D_l$ ($l\geq4$)}
The Dynkin diagram of the simple Lie algebra $D_l\simeq \mathfrak{so}_{2l}^{}$
($l\geq4$) is presented on the picture:
\\
\begin{eqnarray}
\mbox{\begin{picture}(105,40)
\put(-40,0){\circle{5}}
\put(-37,0){\line(1,0){44}}
\put(-43,-15){\fs$\alpha_{1}^{}$}
\put(10,0){\circle{5}}
\put(13,0){\line(1,0){35}}
\put(7,-15){\fs$\alpha_{2}^{}$}
\put(55,0){$\ldots$}
\put(107,0){\line(-1,0){35}}
\put(110,0){\circle{5}}
\put(113,0){\line(1,0){44}}
\put(110,40){\circle{5}}
\put(110,3){\line(0,1){34}}
\put(120,38){\fs$\alpha_l^{}$}
\put(100,-15){\fs$\alpha_{l-2}^{}$}
\put(160,0){\circle{5}}
\put(157,-15){\fs$\alpha_{l-1}^{}$}
\end{picture}}
\nonumber
\end{eqnarray}
\nopagebreak
\\[0pt]
\centerline{\fs Fig. 8. The Dynkin diagram of the Lie algebra
$\mathfrak{so}_{2l}^{}$.}
\\

In the terms of the orthonormalized basis $\epsilon_i^{}$
($i=1,2,\ldots,l$) the root systems $\Pi$, $\Delta_{+}^{}$ and $\Delta$ of
$\mathfrak{so}_{2l}^{}$ are presented as follows:
\begin{eqnarray}
\Pi\!\!&=\!\!&\{\alpha_1^{}\!=\epsilon_1^{}\!-\epsilon_2^{},
\ldots,\epsilon_{l-2}^{}\!-\epsilon_{l-1}^{},\,
\alpha_{l-1}^{}\!=\epsilon_{l-1}^{}\!-\epsilon_{l}^{},\,
\alpha_{l}^{}\!=\epsilon_{l-1}^{}\!+\epsilon_{l}^{}\},
\label{D1}
\\[5pt]
\Delta_{+}\!\!\!\!&=\!\!&\{\epsilon_i^{}\pm\epsilon_j^{}\,|\;1\le i<j\le l\},
\label{D2}
\\[5pt]
\Delta_{}\!\!&=\!\!&\{\pm\epsilon_i^{}\pm\epsilon_j^{}\,|\;i\neq j;
\,i,j=1,2,\ldots,l\}.
\label{D3}
\end{eqnarray}
The root $\theta:=\alpha_1^{}+2\alpha_2^{}+\ldots+2\alpha_{l-2}^{}
+\alpha_{l-1}^{}+\alpha_{l}^{}=\epsilon_{1}^{}+\epsilon_{2}^{}$ is maximal.

The quantum algebra $U_{q}(\mathfrak{so}_{2l}^{})$ is generated by the
Chevalley $U_q(\mathfrak{gl}_l^{})$-elements $e_{i,-i-1}^{}:=
e_{\epsilon_i-\epsilon_{i+1}}^{}$, $e_{i+1,-i}^{}:=e_{\epsilon_{i+1}-
\epsilon_i}$ $(i=1,2,\ldots,l-1)$ with the relations (\ref{A4}) and the
additional elements $e_{l-1,l}^{}:=e_{\epsilon_{l-1}^{}+\epsilon_{l}^{}}$,
$e_{-l,-l+1}^{}:=e_{-\epsilon_{l-1}^{}-\epsilon_{l}^{}}$ with the relations:
\begin{equation}
\begin{array}{rcl}
q^{e_{i,-i}}e_{l-1,l}^{}\,q^{-e_{i,-i}}\!\!&=\!\!&
q^{\delta_{i,l-1}+\delta_{i,l}}e_{l-1,l}^{}~,
\\[5pt]
q^{e_{i,-i}}e_{-l,-l+1}^{}\,q^{-e_{i,-i}}\!\!&=\!\!&
q^{-\delta_{i,l-1}-\delta_{i,l}}e_{-l,-l+1}^{}~,
\\[5pt]
[e_{l-1,l}^{},\,e_{-l,-l+1}^{}]\!\!&=\!\!&
\mbox{\ls$\frac{q^{e_{l-1,-l+1}^{}+e_{l,-l}^{}}\,-
\,q^{-e_{l-1,-l+1}^{}-e_{l,-l}^{}}}
{q\,-\,q^{-1}}$},
\\[5pt]
[e_{i,-j}^{},e_{l-1,l}^{}]\!\!&=\!\!&0\quad\;
(j\neq l\!-\!1,\;|i-j|=1),
\\[5pt]
[e_{i,-j}^{},e_{-l,-l+1}^{}]\!\!&=\!\!&0
\quad(i\neq l\!-\!1,\;|i-j|=1),
\end{array}
\label{D4'}
\end{equation}
\begin{equation}
\begin{array}{rcl}
[[e_{l-2,-l+1}^{},\,e_{l-1,l}^{}]_{q}^{},\,e_{l-1,l}^{}]_{q}^{}
\,\!\!&=\!\!&0~,
\\[5pt]
[[e_{l-1,-l+2}^{},\,e_{-l,-l+1}^{}]_{q}^{},\,e_{-l,-l+1}^{}]_{q}^{}
\!\!&=\!\!&0~,
\\[5pt]
[e_{l-2,-l+1}^{},\,[e_{l-2,-l+1}^{},
\,e_{l-1,l}]_{q}^{}]_{q}^{}\,\!\!&=\!\!&0~,
\\[5pt]
[e_{l-1,-l+2}^{},\,[e_{l-1,-l+2}^{},\,e_{-l+1,-l}]_{q}^{}]_{q}^{}
\,\!\!&=\!\!&0~.
\label{D4}
\end{array}
\end{equation}
The Hopf structure on $U_{q}(\mathfrak{so}_{2l}^{})$ is given by the formulas
(\ref{A5})--(\ref{A7}) and also
\begin{equation}
\begin{array}{rcl}
\Delta_{q}(e_{l-1,l}^{})\!\!&=\!\!&e_{l-1,-l}^{}\otimes1+
q^{e_{l-1,-l+1}^{}+e_{l,-l}^{}}\otimes e_{l-1,l}^{},
\\[5pt]
\Delta_{q}(e_{-l,-l+1}^{})\!\!&=\!\!&e_{-l,-l+1}^{}\otimes
q^{-e_{l-1,-l+1}^{}-e_{l,-l}}+1\otimes e_{-l,-l+1}^{},
\end{array}
\label{D5}
\end{equation}
\begin{equation}
\begin{array}{rcl}
S_{q}(e_{l-1,l}^{})\!\!&=\!\!&-q^{e_{l-1,-l+1}^{}+e_{l,-l}^{}}\,e_{l-1,l}^{},
\\[5pt]
S_{q}(e_{-l,-l+1}^{})\!\!&=\!\!&-e_{-l,-l+1}^{}\,q^{-e_{l-1,-l+1}^{}-e_{l,-l}^{}},
\phantom{aaaaaaaaaaa}
\end{array}
\label{D6}
\end{equation}
\begin{equation}
\varepsilon_{q}(e_{l-1,l}^{})\,=\,\varepsilon_{q}(e_{l+1,-l}^{})\,=\,0~.
\phantom{aaiaaaaaaaaa}
\label{D7}
\end{equation}

For construction of the composite root vectors $e_{i,-j}^{}$
($|i-j|>1$), $e_{i,j}^{}$, $e_{-j,-i}^{}$ ($i<j$, $i\ne l-1$)
we fix the following normal ordering of the positive root system
$\Delta_{+}^{}$: 
\begin{equation}
\begin{array}{llc}
&(\epsilon_1^{}-\epsilon_2^{},\epsilon_1^{}-\epsilon_3^{},\ldots,
\epsilon_1^{}-\epsilon_l^{},\epsilon_1^{}+\epsilon_l^{},
\ldots,\epsilon_1^{}+\epsilon_3^{},\epsilon_{1}^{}+\epsilon_{2}^{}),&
\\[5pt]
&(\epsilon_2^{}-\epsilon_3^{},\epsilon_2^{}-\epsilon_4^{}\ldots,
\epsilon_2^{}-\epsilon_l^{},\epsilon_2^{}+\epsilon_l^{},\ldots,\epsilon_2^{}+
\epsilon_4^{},\epsilon_{2}^{}+\epsilon_{3}^{}),\ldots,&
\\[5pt]
&(\epsilon_{l-2}^{}-\epsilon_{l-1}^{},\epsilon_{l-2}^{}-\epsilon_{l}^{},
\epsilon_{l-2}^{}+\epsilon_{l}^{},\epsilon_{l-2}^{}+\epsilon_{l-1}^{}),
(\epsilon_{l-1}^{}-\epsilon_{l}^{},\epsilon_{l-1}^{}+\epsilon_{l}^{})~.&
\end{array}
\label{D8}
\end{equation}
According to this ordering we set
\begin{equation}
e_{i,-j}^{}\;:=\;[e_{i,-i-1}^{},\,e_{i+1,-j}^{}]_{q^{-1}},
\qquad e_{j,-i}^{}\;:=\;[e_{j,-i-1}^{},\,e_{i+1,-i}^{}]_{q}
\label{D9f}
\end{equation}
for $1\le i\!+\!1<\!j\le l$,
\label{D9b}
\begin{equation}
e_{i,l}^{}\;:=\;[e_{i,-l+1}^{},\,e_{l-1,l}^{}]_{q^{-1}},
\qquad e_{-l,-i}^{}\;:=\;[e_{-l,-l+1}^{},\,e_{l-1,-i}^{}]_{q}
\end{equation}
for $1\le i\leq l\!-\!2$, 
\begin{equation}
e_{i,j}^{}\;:=\;[e_{i,j+1}^{},\,e_{j,-j-l}^{}]_{q^{-1}},
\qquad e_{-j,-i}^{}\;:=\;[e_{j+1,-j}^{},\,e_{-j-1,-i}^{}]_{q}
\label{D9}
\end{equation}
for $1\le i<j\leq l\!-\!2$.

The extended Dynkin diagram of $D_l^{}$ (or the Dynkin diagram of the
corresponding affine Lie algebra $\widehat{\mathfrak{so}}_{2l}^{}$
($l\geq4$)) is presented by the picture: 
\\
\begin{eqnarray}
\mbox{\begin{picture}(130,45)
\put(-50,0){\circle{5}}
\put(-47,0){\line(1,0){44}}
\put(-53,-15){\fs$\alpha_{1}^{}$}
\put(0,0){\circle{5}}
\put(3,0){\line(1,0){44}}
\put(-3,-15){\fs$\alpha_{2}^{}$}
\put(0,40){\circle{5}}
\put(0,3){\line(0,1){34}}
\put(10,38){\fs$\delta-\theta$}
\put(50,0){\circle{5}}
\put(53,0){\line(1,0){35}}
\put(47,-15){\fs$\alpha_{3}^{}$}
\put(93,0){$\ldots$}
\put(147,0){\line(-1,0){35}}
\put(150,0){\circle{5}}
\put(153,0){\line(1,0){44}}
\put(150,40){\circle{5}}
\put(150,3){\line(0,1){34}}
\put(160,38){\fs$\alpha_{l}^{}$}
\put(140,-15){\fs$\alpha_{l-2}^{}$}
\put(200,0){\circle{5}}
\put(197,-15){\fs$\alpha_{l-1}^{}$}
\end{picture}}
\label{D10'}
\nonumber
\end{eqnarray}
\nopagebreak
\\[0pt]
\centerline{\fs Fig. 9. The Dynkin diagram of the Lie algebra
$\widehat{\mathfrak{so}}_{2l}^{}$.}
\\

Specializing the general formulas (\ref{tolsD5})--(\ref{tolsD8}) to
the case $\mathfrak{g}=\mathfrak{so}_{2l}^{}$ with the vector
\begin{equation}
\tilde{e}_{-\theta}^{}=q^{e_{1,-1}^{}\!+e_{2,-2}^{}}e_{-2,-1}~,
\label{D10}
\end{equation}
after some calculations we obtain the following result.
\begin{theorem}
The Drinfeldian $D_{q\eta}(\mathfrak{so}_{2l}^{})$ ($l\ge 4$) is generated
(as a unital associative algebra over $\mathbb C[[\log q,\eta]]$) by the
algebra $U_{q}(\mathfrak{so}_{2l}^{})$ and the elements
$\xi_{\delta-\theta}^{}$, $q^{\pm\hat{c}}:=q^{\pm h_\delta}$ with the
relations:
\begin{eqnarray}
[q^{\pm\hat{c}}\!,\,{\rm everything}]\!\!&=\!\!&0~,
\label{D11}
\\[4pt]
q^{\pm e_{i,-i}}\xi_{\delta-\theta}^{}\!\!&=\!\!&q^{\mp1}
\xi_{\delta-\theta}^{}q^{\pm e_{i,-i}}\quad {\rm for}\;\, i=1,2,
\label{D12}
\\[4pt]
q^{\pm e_{i,-i}}\xi_{\delta-\theta}^{}\!\!&=\!\!&\xi_{\delta-\theta}^{}\,
q^{\pm e_{i,-i}}\qquad\;{\rm for}\;\, i=3,4,\ldots,l,
\label{D13}
\\[3pt]
[e_{2,-1}^{},\,\xi_{\delta-\theta}^{}]\!\!&=\!\!&\eta\,
q^{e_{1,-1}^{}+e_{2,-2}^{}}\sum_{k=3}^{l}(-1)^kq^{k-3}e_{-k,-1}^{}e_{k,-1}^{},
\label{D14}
\\[3pt]
[e_{1,-2}^{},\,\xi_{\delta-\theta}^{}]\!\!&=\!\!&
\eta\,q^{2e_{2,-2}^{}}\sum_{k=3}^{l}(-1)^kq^{k-3}e_{-k,-2}^{}e_{k,-2}^{},
\label{D15}
\end{eqnarray}
\begin{eqnarray}
[e_{3,-2}^{},\,\xi_{\delta-\theta}^{}]\!\!&=\!\!&0~,
\\[5pt]
\label{D16}
[e_{2,-3}^{},\,[e_{2,-3}^{},\,\xi_{\delta-\theta}^{}]_{q}]_{q}\!\!&=\!\!&0~,
\label{D17}
\\[5pt]
[[e_{2,-3}^{},\,\xi_{\delta-\theta}^{}]_{q},\,
\xi_{\delta-\theta}^{}]_{q}\!\!&=\!\!&0~,
\label{D18}
\\[5pt]
[e_{i,-j}^{},\,\xi_{\delta-\theta}^{}]\!\!&=\!\!&0
\quad\; 
(3\leq i,j\leq l,\;\,\mid i-j\mid=1),
\\[5pt]
[e_{l-1,l}^{}\,\xi_{\delta-\theta}^{}]\!\!&=\!\!&0~,
\\[5pt]
[e_{-l,-l+1}^{}\,\xi_{\delta-\theta}^{}]\!\!&=\!\!&0~.
\label{D19}
\end{eqnarray}
\end{theorem}
By setting $q=1$ in (\ref{D11})--(\ref{D19}), we obtain the defining
relations of the Yangian $Y_{\eta}(\mathfrak{so}_{2l}^{})$
in the Chevalley basis. This result is formulated as the theorem.
\begin{theorem}
The Yangian $Y_{\eta}(\mathfrak{so}_{2l}^{})$ ($l\geq4$) is generated
(as a unital associative algebra over $\mathbb C[\eta]$) by the algebra
$U(\mathfrak{so}_{2l}^{})$ and the elements $\xi_{\delta-\theta}^{}$,
$\hat{c}$ with the relations:
\begin{eqnarray}
[\hat{c},\,{\rm everything}]\!\!&=\!\!&0~,
\label{D20}
\\[4pt]
[e_{i,-i},\,\xi_{\delta-\theta}^{}]\!\!&=\!\!&\xi_{\delta-\theta}^{}
\quad {\rm for}\;\,i=1,2,
\label{D21}
\\[4pt]
[e_{i,-i}^{}\xi_{\delta-\theta}^{}]\!\!&=\!\!&0
\qquad\;\,{\rm for}\;\, i=3,4,\ldots,l,
\label{D22}
\\[3pt]
[e_{2,-1}^{},\,\xi_{\delta-\theta}^{}]\!\!&=\!\!&
\eta\sum_{k=3}^{l}(-1)^ke_{-k,-1}^{}e_{k,-1}^{},
\label{D23}
\\[3pt]
[e_{1,-2}^{},\,\xi_{\delta-\theta}^{}]\!\!&=\!\!&
\eta\sum_{k=3}^{l}(-1)^ke_{-k,-2}^{}e_{k,-2}^{},
\label{D24}
\\[3pt]
[e_{3,-2}^{},\,\xi_{\delta-\theta}^{}]\!\!&=\!\!&0~,
\label{D25}
\\[5pt]
[e_{2,-3}^{},\,[e_{2,-3}^{},\,\xi_{\delta-\theta}^{}]]\!\!&=\!\!&0~,
\label{D26}
\\[5pt]
[[e_{2,-3}^{},\,\xi_{\delta-\theta}^{}],\,\xi_{\delta-\theta}^{}]\!\!&=\!\!&0~,
\label{D27}
\\[5pt]
[e_{i,-j}^{},\,\xi_{\delta-\theta}^{}]\!\!&=\!\!&0
\quad\; 
(3\leq i,j\leq l,\;\,|i-j|=1),
\label{D28'}
\end{eqnarray}
\begin{eqnarray}
[e_{l-1,l}^{}\,\xi_{\delta-\theta}^{}]\!\!&=\!\!&0~,
\\[4pt]
[e_{-l,-l+1}^{}\,\xi_{\delta-\theta}^{}]\!\!&=\!\!&0~.
\label{D28}
\end{eqnarray}
\end{theorem}

\bibliographystyle{amsalpha}

\begin{thebibliography}{A}

\bibitem{BD}
A.A. Belavin and V.G. Drinfeld,
\textit{On solutions of classical Yang-Baxter equation for simple Lie algebras},
Funct. Analysis and its Appl. \textbf{16}, No. 3 (1982), 1--29.

\bibitem{D1}
V.G. Drinfeld,
\textit{Hopf algebras and quantum Yang-Baxter equation},
Soviet Math. Dokl. \textbf{283} (1985), 1060--1064.

\bibitem{D2}
V.G. Drinfeld,
\textit{A new realization of Yangians and quantized affine algebras},
Soviet Math. Dokl. \textbf{32} (1988), 212--216.

\bibitem{D3}
V.G. Drinfeld,
\textit{Quantum groups}, Proc. ICM-86 (Berkeley USA), vol. 1,
Amer. Math. Soc. Providence, RI, 1987, pp. 798--820.

\bibitem{Kac}
V.G. Kac,
\textit{Infinite dimensional Lie algebras},
Cambridge University Press, Cambrige, USA, 1985.

\bibitem{KT1}
S.M. Khoroshkin and V.N. Tolstoy,
\textit{Universal R-matrix for quantized (super)algebras},
Commun. Math. Phys. \textbf{141}, No. 3 (1991), 599--617.

\bibitem{KT2}
S.M. Khoroshkin and V.N. Tolstoy,
\textit{Twisting of quantum (super)algebras. Connection of Drinfeld's
and Cartan-Weyl realizations for quantum affine algebras},
MPIM preprint, MPI/94-23, pp. 1--29 (Bonn, 1994);
{\tt arXiv:hep-th/9404036}.

\bibitem{KT3}
S.M. Khoroshkin and V.N. Tolstoy,
\textit{Yangian double},
Lett. Math. Phys. \textbf{36} (1996), 373--402:
{\tt arXiv:hep-th/9406194}.

\bibitem{MNO}
A. Molev, M. Nazarov and G. Olshanski,
\textit{Yangians and classical Lie algebras},
Uspekhi Math. Nauk \textbf{51} (1996), 27--104.

\bibitem{RTF}
N.Yu. Reshetikhin, L.A. Takhtadjan and L.D. Faddeev,
\textit{Quantization of Lie groups and Lie algebras},
Lelingrad Math. J. \textbf{1} (1990), 193--225.

\bibitem{T1}
V.N. Tolstoy,
\textit{Extremal projectors for quantized Kac-Moody superalgebras
and some of their applications},
Lecture Notes in Phys., (Springer, Berlin), \textbf{370}
(1990), 118--125.

\bibitem{T2}
V.N. Tolstoy,
\textit{Connection between Yangians and Quantum Affine Algebras},
Proc. of the X-th Max Born Symposium,
(Wroclav, 1996), eds. J. Lukierski, M. Mozrzymas,
PWN - Polish Sci. Publishers, Warszawa, 1997, pp. 99--117.

\bibitem{T3}
V.N. Tolstoy,
\textit{Drinfeldians}, Proc. of an International Workshop "Lie theory
and its application in physics II",
(Clausthal, 1997),
eds. H.-D. Doebner, V.K. Dobrev and J. Hilgert,
World Sci. Publishing, River Edge, NJ, 1998,
pp. 225--337;
{\tt math.QA/9803008}.

\bibitem{T4}
V.N. Tolstoy,
\textit{From Quantum Affine Superalgebras to super-Drinfeldians
and super-Yangians}, Proc. of JINR Workshop
"Supersymmetries and Quantum Symmetries", (Dubna, 1999),
eds. E. Ivanov, S. Krivonos,
Joint Institute for Nuclear Research, Dubna, 2000,
pp. 431--439.

\bibitem{T5}
V.N. Tolstoy,
\textit{Super-Drinfeldians and super-Yangians of Lie superalgebra
of type $A(n|m)$},
Physics of Atomic Nuclei \textbf{64}, No. 12 (2001), 2179--2184.
\end{thebibliography}

\end{document}